\documentclass[conference]{IEEEtran}
\IEEEoverridecommandlockouts
\usepackage{cite}
\usepackage{amsmath,amssymb,amsfonts}
\usepackage{algorithmic}
\usepackage{subcaption}
\usepackage{graphicx}
\usepackage{textcomp}
\usepackage[dvipsnames]{xcolor}
\usepackage{subcaption}
\usepackage[linesnumbered,ruled,vlined]{algorithm2e}
\def\BibTeX{{\rm B\kern-.05em{\sc i\kern-.025em b}\kern-.08em
    T\kern-.1667em\lower.7ex\hbox{E}\kern-.125emX}}
\definecolor{Yellow2}{HTML}{C6BD2F}

\newcommand{\rev}[1]{{\color{black}{#1}}}
\newcommand{\redtext}[1]{{\color{red}{#1}}}
\newcommand{\greentext}[1]{{\color{PineGreen}{#1}}}
\newcommand{\bluetext}[1]{{\color{RoyalBlue}{#1}}}
\newcommand{\yellowtext}[1]{{\color{Yellow2}{#1}}}
\newcommand{\orangetext}[1]{{\color{orange}{#1}}}

\newcommand\tp[2][-6]{{#2}^{\mkern#1mu\top}}

\begin{document}

\title{
Scalable Linear Time Dense Direct Solver for 3-D Problems Without Trailing Sub-Matrix Dependencies
\thanks{This work was supported by JSPS KAKENHI Grant Number JP20K20624,JP21H03447.This work is supported by ”Joint Usage/Research Center for Interdisciplinary Large-scale Information Infrastructures” in Japan (Project ID: jh210024-NAHI)}
}

\author{\IEEEauthorblockN{Qianxiang Ma}
\IEEEauthorblockA{\textit{School of Computing} \\
\textit{Tokyo Institute of Technology}\\
Tokyo, Japan \\
ma@rio.gsic.titech.ac.jp}
\and
\IEEEauthorblockN{Sameer Deshmukh}
\IEEEauthorblockA{\textit{School of Computing} \\
\textit{Tokyo Institute of Technology}\\
Tokyo, Japan \\
sameer.deshmukh@rio.gsic.titech.ac.jp}
\and
\IEEEauthorblockN{Rio Yokota}
\IEEEauthorblockA{\textit{Global Scientific Information and Computing Center} \\
\textit{Tokyo Institute of Technology}\\
Tokyo, Japan \\
rioyokota@gsic.titech.ac.jp}
}

\maketitle

\begin{abstract}
\rev{Factorization of large dense matrices are ubiquitous in engineering and data science applications, \textit{e.g.} preconditioners for iterative boundary integral solvers, frontal matrices in sparse multifrontal solvers, and computing the determinant of covariance matrices.}
HSS and $\mathcal{H}^2$-matrices are hierarchical low-rank matrix formats that can reduce the complexity of factorizing \rev{such} dense matrices from $\mathcal{O}(N^3)$ to $\mathcal{O}(N)$. 
For HSS matrices, it is possible to remove the dependency on the \rev{trailing matrices} during Cholesky/LU factorization, which results in a highly parallel algorithm. 
However, the weak admissibility of HSS \rev{causes the rank of off-diagonal blocks to grow for 3-D problems, and the method is no longer $\mathcal{O}(N)$}.
On the other hand, the strong admissibility of $\mathcal{H}^2$-matrices allows it to handle 3-D problems \rev{in $\mathcal{O}(N)$}, but introduces a dependency on the trailing matrices. 
In the present work, we \rev{pre-compute the fill-ins and integrate them into the shared basis}, which allows us to remove the dependency on trailing-matrices even for $\mathcal{H}^2$-matrices.
\rev{Comparisons with a block low-rank factorization code LORAPO showed a maximum speed up of 4,700x for a 3-D problem with complex geometry.}
\end{abstract}


\section{Introduction}
Factorization of dense matrices has been at the heart of high performance computing where the high performance LINPACK (HPL) benchmark has been used to track the performance of the Top500 supercomputers for nearly 30 years.
Even though HPL is often criticized for not being representative of the actual workloads in HPC, there exist many applications that require the factorization of large dense matrices.
For example, preconditioners for iterative boundary integral solvers~\cite{takahashiParallelizationInverseFast2020}, \rev{frontal matrices} in sparse multifrontal solvers~\cite{liuSparseApproximateMultifrontal2020,amestoyPerformanceScalabilityBlock2019}, and computing the determinant of \rev{covariance} matrices in statistics~\cite{litvinenkoLikelihoodApproximationHierarchical2019}.
The dense matrices that arise in these problems have a low-rank structure that can be exploited to perform approximate factorization in linear time.
The accuracy of the approximation is controllable by adjusting the rank or the admissibility.
The admissibility condition determines how far to subdivide the blocks, and which ones to consider as dense blocks.
The hierarchical semi-separable (HSS) matrix~\cite{chandrasekaranFastSolverHSS2006} has weak admissibility, where only the diagonal blocks are subdivided recursively, and off-diagonal blocks are approximated with low-rank matrices.
\rev{Since the numerical rank of each block is determined by the proximity of the sub-domains in the underlying geometry, for 3-D problems this large off-diagonal block in HSS will contain many sub-blocks that have large rank.
This causes the rank of off-diagonal blocks to grow with the problems size, and the method would no longer have linear complexity.}
The $\mathcal{H}^2$-matrix~\cite{hackbuschMatrices2000} on the other hand has strong admissibility, where the off-diagonal blocks can be subdivided \rev{further so that each block can maintain a somewhat constant rank that is independent of the problem size.
For the off-diagonal sub-blocks that have close proximity in the geometry will be subdivided until the leaf-level is reached, at which point it will be treated as a dense block.
There are only a constant number of such off-diagonal dense blocks per row/column, due to the number of neighboring boxes in the geometry being constant and independent of the problem size.
This allows the $\mathcal{H}^2$-matrix to achieve linear complexity even for 3-D problems.}
Both, HSS and $\mathcal{H}^2$-matrices share the basis among the low-rank blocks in the same block row/column.
These shared bases are also nested between the levels so that the bases of large low-rank blocks do not need to be stored explicitly.
Variants that use independent bases for each low-rank block also exist, \textit{e.g.} HODLR (weak admissibility)~\cite{ambikasaranFastAlgorithmsDense2013} and $\mathcal{H}$-matrix (strong admissibility)~\cite{hackbuschSparseMatrixArithmetic1999}. 
There are also non-hierarchical variants such as BLR (independent basis)~\cite{amestoyImprovingMultifrontalMethods2015} and $\mathrm{BLR}^2$ (shared basis)~\cite{ashcraftBlockLowRankMatrices2021}.
Table \ref{tab:structures} summarizes the various types of structures categorized by their basis and admissibility, along with the factorization complexity.

\begin{table}[b]
    \centering
    \begin{tabular}{c|c|c|c}
        Low-rank structure & Basis & Admissibility & Complexity\\
        \hline
        BLR~\cite{amestoyImprovingMultifrontalMethods2015} & Independent & Strong or weak & $O(N^2)$ \\
        $\mathrm{BLR}^2$~\cite{ashcraftBlockLowRankMatrices2021} & Shared & Strong or weak & $O(N^{1.8})$ \\
        HODLR~\cite{ambikasaranFastAlgorithmsDense2013} & Independent & Weak & $O(N\log^2 N)$ \\
        $\mathcal{H}$-matrix~\cite{hackbuschSparseMatrixArithmetic1999} & Independent & Strong & $O(N\log^2 N)$ \\
        HSS~\cite{chandrasekaranFastSolverHSS2006} & Shared & Weak & $O(N)$ \\
        $\mathcal{H}^2$-matrix~\cite{hackbuschMatrices2000} & Shared & Strong & $O(N)$ \\
    \end{tabular}
    \caption{List of Different Low-rank Structures}
    \label{tab:structures}
\end{table}

The different structures each have their pros and cons. 
The non-hierarchical structures such as BLR and $\mathrm{BLR}^2$ are obviously much simpler to implement, but cannot achieve the near-linear complexity that hierarchical structures enjoy.
It is worth noting that \rev{frontal matrices} in multifrontal solvers have a size of \rev{$\mathcal{O}(N^{\frac{2}{3}})$} \rev{for 3-D problems}, so reducing the complexity of factorizing this part to $\mathcal{O}(N^2)$ \rev{with BLR} is enough to achieve an overall complexity of \rev{$\mathcal{O}(N^{\frac{4}{3}})$} for the multifrontal solver.
Sharing the basis among the block rows/columns reduces the memory consumption, but requires the basis for the whole block row/column to be updated every time a trailing sub-block is updated.
The ULV factorization for HSS matrices~\cite{chandrasekaran_ulv_2006} avoids this problem by factorizing only the redundant part of the dense blocks at each level.

\begin{figure}[t]
  \centering
  \includegraphics[width=\linewidth]{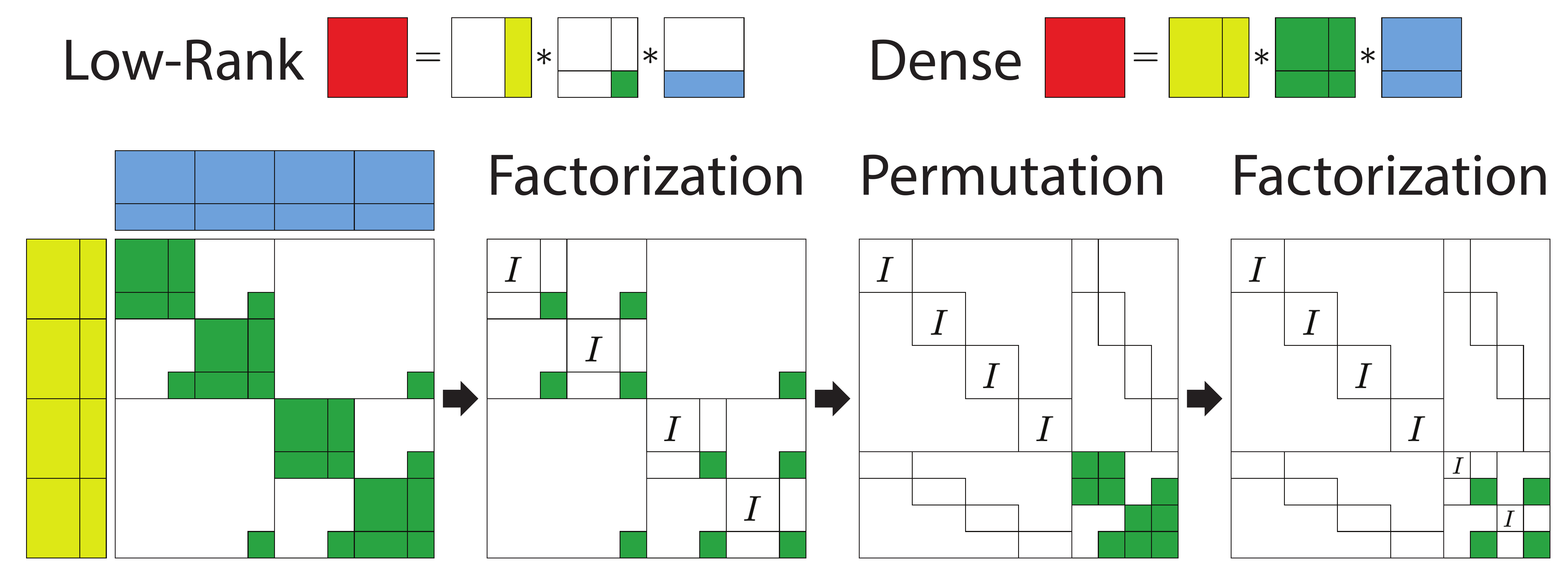}
  \caption{$\mathcal{O}(N)$ Dense Factorization Without Trailing Sub-Matrix Dependencies}
  \label{fig:HSS-ULV}
\end{figure}

As shown at the top of Fig. \ref{fig:HSS-ULV}, low-rank blocks can be decomposed into the \yellowtext{shared column basis}, \greentext{skeleton matrix}, and  \bluetext{shared row basis}.
We can then add back the redundant part of the basis to the \yellowtext{shared column basis} and \bluetext{shared row basis} to form a square matrix.
These bases are used to compute the \greentext{skeleton matrix} for the dense block, which are partitioned into four parts.
This allows us to decompose a dense matrix into a sparsifed structure shown on the leftmost side of Fig. \ref{fig:HSS-ULV}.
This block sparse matrix can be factorized without trailing sub-matrix dependencies, if we leave the skeleton matrices to be factorized at the next level.
This has huge implications with regards to both the complexity and parallelism of the algorithm.
The sparsification results in $\mathcal{O}(N)$ compleixity, and the removal of dependencies \rev{allows the diagonal blocks to be factorized in parallel}.
Concepts such as right-looking and left-looking become irrelevant, since there are no GEMM operations nor trailing sub-matrix dependencies required in this method.
Runtime systems such as StarPU~\cite{augonnetStarPUUnifiedPlatform2011} and PaRSEC~\cite{bosilcaDAGuEGenericDistributed2012} and coloring schemes to extract parallelism are also unnecessary since the method has no dependencies.

However, the procedure illustrated in Fig. \ref{fig:HSS-ULV} only works for HSS matrices, which has \rev{suboptimal complexity for} problems with 3-D geometry, as mentioned earlier.
On the other hand, $\mathcal{H}^2$-matrices \rev{can still achieve linear complexity even for 3-D problems.}
However, having dense blocks in the off-diagonal will result in fill-in during the factorization, which \rev{introduces} the dependency on the trailing sub-matrices to this otherwise highly parallel factorization method.
In the present work, we \rev{pre-compute the fill-ins and integrate them into the shared basis}, in order to factorize an $\mathcal{H}^2$-matrix without any dependency on the trailing sub-matrices.
Unlike existing methods that rely on coloring to extract parallelism from an $\mathcal{H}^2$-matrix factorization, our method is inherently parallel.

We make the following contributions in the present work:
\begin{itemize}
    \item We extend the HSS-ULV factorization to $\mathcal{H}^2$-matrices, which \rev{has $\mathcal{O}(N)$ complexity even} for 3-D problems.
    \item We do this while retaining the parallelism of HSS-ULV by \rev{pre-computing the fill-ins and including them in the shared basis}.
    \item \rev{This is the first method that can factorize dense matrices arising from 3-D problems in linear time without trailing sub-matrix dependencies.}
\end{itemize}

\section{Structured Low-Rank Factorization}\label{sec:structured}

\begin{figure}[t]
  \centering
  \includegraphics[width=\linewidth]{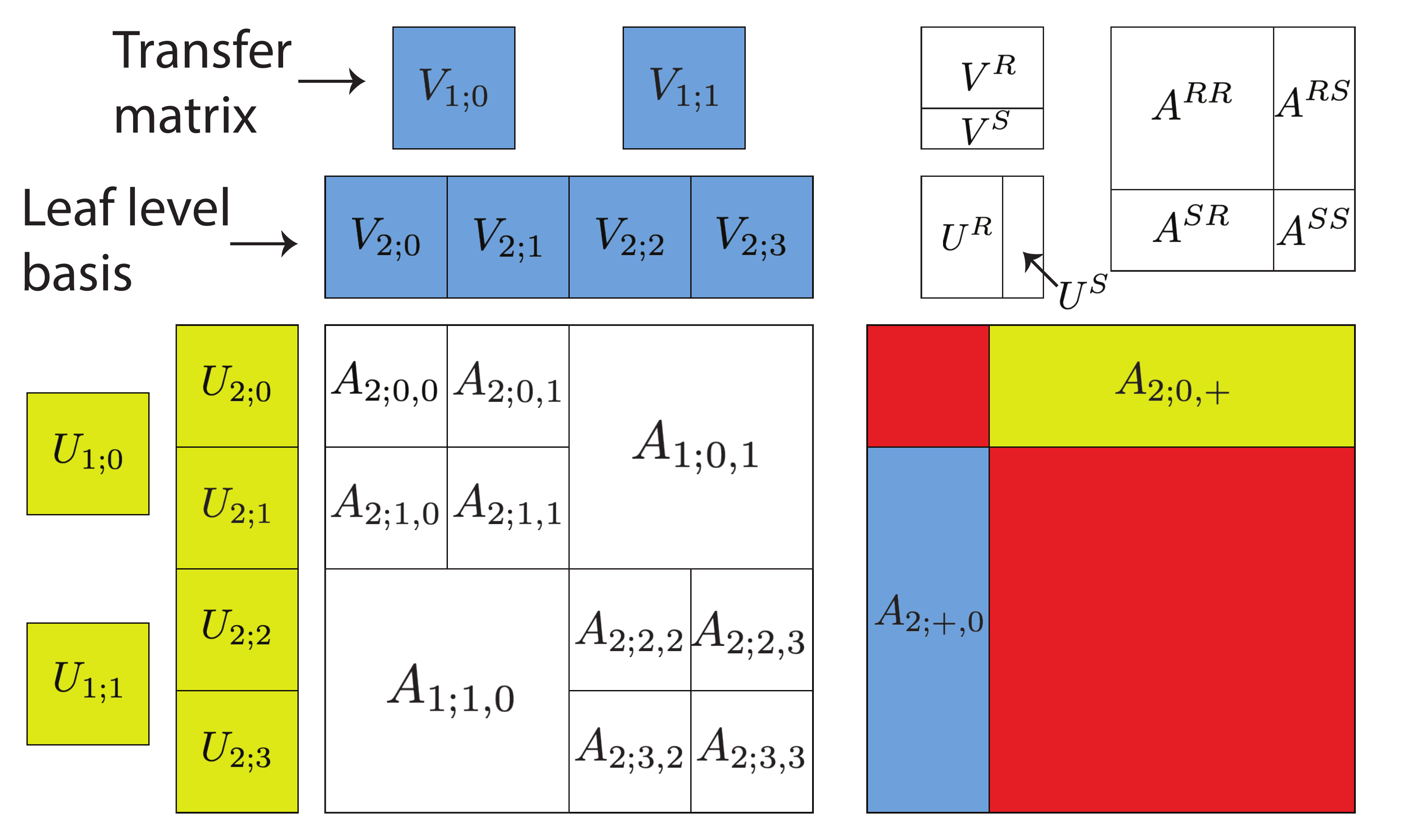}
  \caption{Index notation for hierarchical matrices}
  \label{fig:notation}
\end{figure}

\subsection{Notation}
We first define the common notations we use throughout this paper in Fig. \ref{fig:notation}.
A sub-block of a hierarchically sub-divided matrix $A$ is indexed using $A_{level;row,column}$.
The shared column basis and row basis are indexed by $U_{level;row}$ and $V_{level;column}$, respectively.
At the leaf level this notation is used for the shared basis and at other levels it is used to denote the transfer matrix.
For example, matrix $A_{1;0,1}$ can be written in the following hierarchical low-rank form
\begin{equation}
\redtext{A_{1;0,1}}=
\yellowtext{
\begin{bmatrix}
U_{2;0} & 0 \\
0 & U_{2;1}
\end{bmatrix}
U_{1;0}}
\greentext{S_{1;0,1}}
\bluetext{V_{1;1}
\begin{bmatrix}
V_{2;2} & 0 \\
0 & V_{2;3}
\end{bmatrix}}.
\label{eq:transfer}
\end{equation}
where \greentext{$S_{1;0,1}$} is the \greentext{skeleton matrix} shown inside the corresponding large off-diagonal block in Fig. \ref{fig:HSS-ULV}. We also introduce a convenient notation for expressing the concatenation of all low-rank blocks in a given row/column as shown on the right in Fig. \ref{fig:notation}. This is used in the first step of the algorithm to compute the shared bases
\begin{align}
\yellowtext{[U^S_{2;0}\ U^R_{2;0}]},R &= \mathrm{QR}(\yellowtext{A_{2;0,+}}) \\
\bluetext{[\tp{V_{2;0}^S}\ \tp{V_{2;0}^R}]},R &= \mathrm{QR}(\bluetext{\tp{A_{2;+,0}}}),
\end{align}
where the $S$ and $R$ superscripts represent the skeleton part and redundant part of the basis, respectively.
The QR() represents a column pivoted or rank revealing QR, where the $R$ matrix is not used.
This can be replaced by an interpolative decomposition~\cite{martinssonRandomizedAlgorithmDecomposition2011} if that is preferred,
\rev{but the loss of orthogonality will lead to other complications in this case}.
In order for the ULV factorization to work, we need to permute the skeleton part and redundant part to look like the shapes shown in Figs. \ref{fig:HSS-ULV} and \ref{fig:notation}.
Using this split in the column and row basis, any dense or low-rank matrix can be subdivided into $A^{RR}$, $A^{RS}$, $A^{SR}$, and $A^{SS}$, as shown in the upper right corner of Fig. \ref{fig:notation}.
Using these notations, the dense matrix $A$ can be decomposed as
\begin{equation}
\redtext{A}=\yellowtext{[U^R\ U^S]}
\greentext{
\begin{bmatrix}
S^{RR} & S^{RS} \\
S^{SR} & S^{SS}
\end{bmatrix}}
\bluetext{
\begin{bmatrix}
V^R \\
V^S
\end{bmatrix}}.
\end{equation}

This is what is shown \rev{in the top-right} of Fig. \ref{fig:BLR2-ULV}.
Similarly, a low-rank matrix $A$ can be decomposed as
\begin{align}
\redtext{A}
&=\yellowtext{[U^R\ U^S]}
\greentext{
\begin{bmatrix}
0 & 0 \\
0 & S^{SS}
\end{bmatrix}}
\bluetext{
\begin{bmatrix}
V^R \\
V^S
\end{bmatrix}}.
\end{align}
The bases $\yellowtext{[U^R\ U^S]}$ and $\bluetext{[V^R\ V^S]}$ are shared among both the dense and low-rank blocks in that entire row and column, respectively.
This allows us to perform a $\yellowtext{U}\greentext{S}\bluetext{V}$ decomposition of the entire matrix, which is what is shown in the leftmost figure in Figs. \ref{fig:HSS-ULV} \rev{and \ref{fig:BLR2-ULV}}.
With these notations, we can proceed to describe our proposed method along with some existing ones.

\subsection{BLR$^2$-ULV factorization}

\begin{figure}[t]
  \centering
  \includegraphics[width=\linewidth]{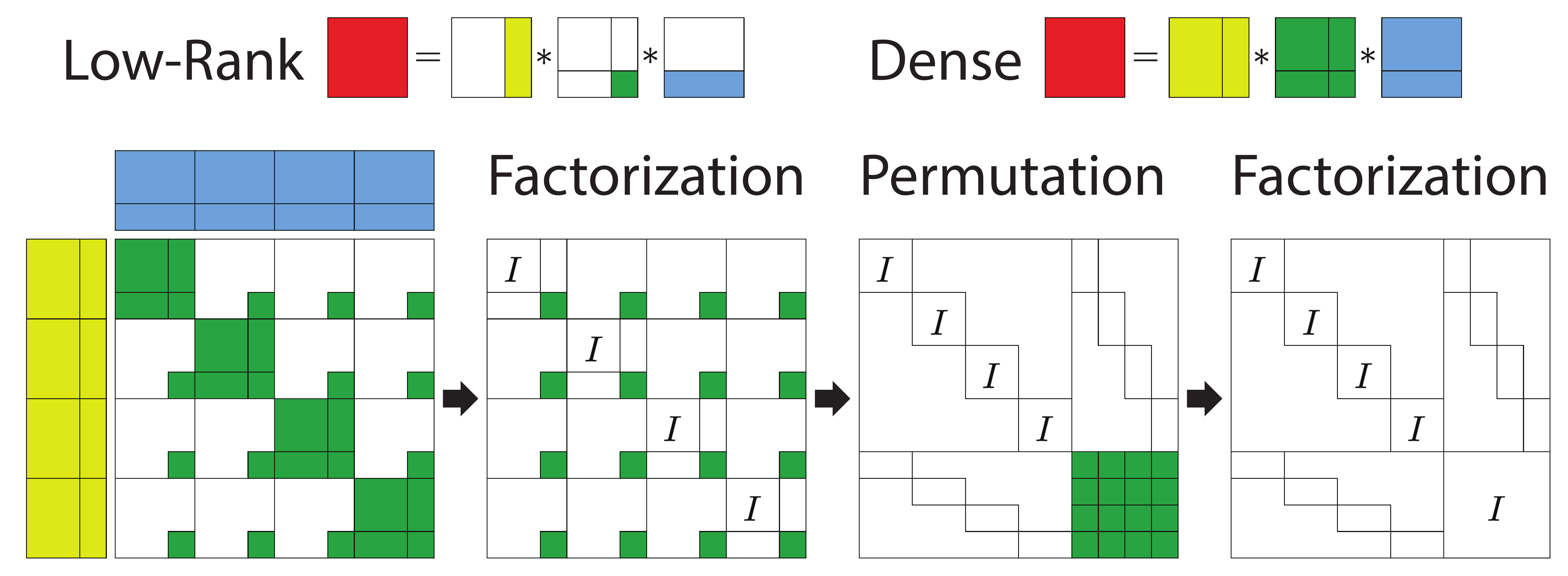}
  \caption{Flow of the BLR$^2$-ULV factorization.}
  \label{fig:BLR2-ULV}
\end{figure}

We will start out with a description of a ULV factorization for the BLR$^2$ structure, which \rev{a non-hierarchical version of the HSS}.
In the following sections, we will add hierarchy to form the HSS-ULV factorization, and then extend that to strong admissibility to obtain the $\mathcal{H}^2$-ULV factorization that we use in our current work.

Similar to the HSS-ULV factorization in Fig. \ref{fig:HSS-ULV} the BLR$^2$-ULV factorization is described in Fig. \ref{fig:BLR2-ULV}. We first obtain the shared bases by performing
\begin{align}
\yellowtext{[U^S_{2;i}\ U^R_{2;i}]},R &= \mathrm{QR}(\yellowtext{A_{2;i,+}}) \label{eq:u}\\ 
\bluetext{[\tp{V_{2;j}^S}\ \tp{V_{2;j}^R}]},R &= \mathrm{QR}(\bluetext{\tp{A_{2;+,j}}}), \label{eq:v}
\end{align}
for all rows of $i$ and all columns of $j$. We then compute all the \greentext{$S$} matrices of the dense diagonal blocks as
\begin{equation}
\greentext{
\begin{bmatrix}
S^{RR}_{2;i,i} & S^{RS}_{2;i,i} \\
S^{SR}_{2;i,i} & S^{SS}_{2;i,i}
\end{bmatrix}}
=\yellowtext{\tp{[U^R_{2;i}\ U^S_{2;i}]}}
\redtext{A_{2;i,i}}
\bluetext{
\tp{
\begin{bmatrix}
V^R_{2;i} \\
V^S_{2;i}
\end{bmatrix}}},
\label{eq:rs}
\end{equation}
for all the diagonal blocks \redtext{$A_{2;,i,i}$}. For the low-rank blocks computing the \greentext{$S$} matrices is simply
\begin{equation}
\greentext{S^{SS}_{2;i,j}}
=\yellowtext{\tp{U^S_{2;i}}}
\redtext{A_{2;i,j}}
\bluetext{
\tp{
V^S_{2;j}}},
\label{eq:ss}
\end{equation}
For all the off-diagonal blocks \redtext{$A_{2;,i,j}$}, where $i\ne j$.
This allows us to decompose for example, the sub-matrix \redtext{$A_{1;0,0}$} into
\begin{equation}
\redtext{A_{1;0,0}}=
\yellowtext{
\setlength\arraycolsep{0pt}
\begin{bmatrix}
U_{2;0} & 0\\
0 & U_{2;1}
\end{bmatrix}}
\greentext{
\begin{bmatrix}
S_{2;0,0} & S_{2;0,1}\\
S_{2;1,0} & S_{2;1,1}
\end{bmatrix}}
\bluetext{
\setlength\arraycolsep{0pt}
\begin{bmatrix}
V_{2;0} & 0\\
0 & V_{2;1}
\end{bmatrix}}.
\label{eq:USV}
\end{equation}
Note that although the \yellowtext{column basis} and \bluetext{row basis} seem like tall skinny matrices in Figs. \ref{fig:HSS-ULV} and \ref{fig:BLR2-ULV}, they are actually block diagonal matrices as shown in Eq. (\ref{eq:USV}).
The \yellowtext{column} and \bluetext{row} basis remain constant throughout the factorization. When a Cholesky factorization is performed only on the \greentext{S} matrix of the $\yellowtext{U}\greentext{S}\bluetext{V}$ decomposition it is called the ULV factorization. 
Actually, calling it a \yellowtext{U}\greentext{LL$^\top$}\bluetext{V} factorization might be more descriptive.
Although the original ULV factorization uses a Cholesky factorization, we may extend this to the LU factorization as well, which is what we will show in the following description.

The LU factorization on the \greentext{S} matrix is done by first eliminating the $S^{RR}_{2;i,i}$ blocks on the diagonal
\begin{equation}
\hat{L}^{RR}_{2;i,i},\hat{U}^{RR}_{2;i,i} = \mathrm{LU}(S^{RR}_{2;i,i}).
\label{eq:lu}
\end{equation}
We use $\hat{U}$ to distinguish the upper triangular matrix from the column basis \yellowtext{$U$}.
This is followed by an elimination of the $S^{RS}_{2;i,i}$ and $S^{SR}_{2;i,i}$ blocks.
\begin{align}
\hat{U}^{RS}_{2;i,i} &= (\hat{L}^{RR}_{2;i,i})^{-1}S^{RS}_{2;i,i}\\
\hat{L}^{SR}_{2;i,i} &= S^{SR}_{2;i,i}(\hat{U}^{RR}_{2;i,i})^{-1}.
\end{align}
Finally, the elimination for the block is completed by computing
\begin{equation}
S^{SS}_{2;i,i} = S^{SS}_{2;i,i} -  \hat{L}^{SR}_{2;i,i}\hat{U}^{RS}_{2;i,i}.
\end{equation}
The elimination of these blocks for different $i$ can be done in parallel since there are no dependencies among them.
Following the factorization phase, the $S^{SS}_{2;i,j}$ blocks that remain to be factorized are clustered in the bottom-left corner, as shown in the permutation phase in Fig. \ref{fig:BLR2-ULV}.
For the BLR$^2$-ULV factorization, this entire remaining block is eliminated as a single dense matrix, which completes the LU factorization.
\begin{equation}
\hat{L}^{SS}_{2;1:4,1:4},\hat{U}^{SS}_{2;1:4,1:4}=
LU\left(
\setlength\arraycolsep{0pt}
\begin{bmatrix}
S^{SS}_{2;0,0} & S^{SS}_{2;0,1} & S^{SS}_{2;0,2} & S^{SS}_{2;0,3}\\
S^{SS}_{2;1,0} & S^{SS}_{2;1,1} & S^{SS}_{2;1,2} & S^{SS}_{2;1,3}\\
S^{SS}_{2;2,0} & S^{SS}_{2;2,1} & S^{SS}_{2;2,2} & S^{SS}_{2;2,3}\\
S^{SS}_{2;3,0} & S^{SS}_{2;3,1} & S^{SS}_{2;3,2} & S^{SS}_{2;3,3}
\end{bmatrix}
\right).
\label{eq:BLR2final}
\end{equation}
The forward elimination and backward substitution of the resulting
\begin{equation}
\yellowtext{U}\greentext{\hat{L}\hat{U}}\bluetext{V}x = b
\end{equation}
can be performed in the following three steps
\begin{align}
\greentext{\hat{L}}z &= \yellowtext{\tp{U}}b\\
\greentext{\hat{U}}y &= z\\
x &= \bluetext{\tp{V}}y.
\end{align}

\subsection{HSS-ULV factorization}
\begin{figure}[t]
  \centering
  \includegraphics[width=\linewidth]{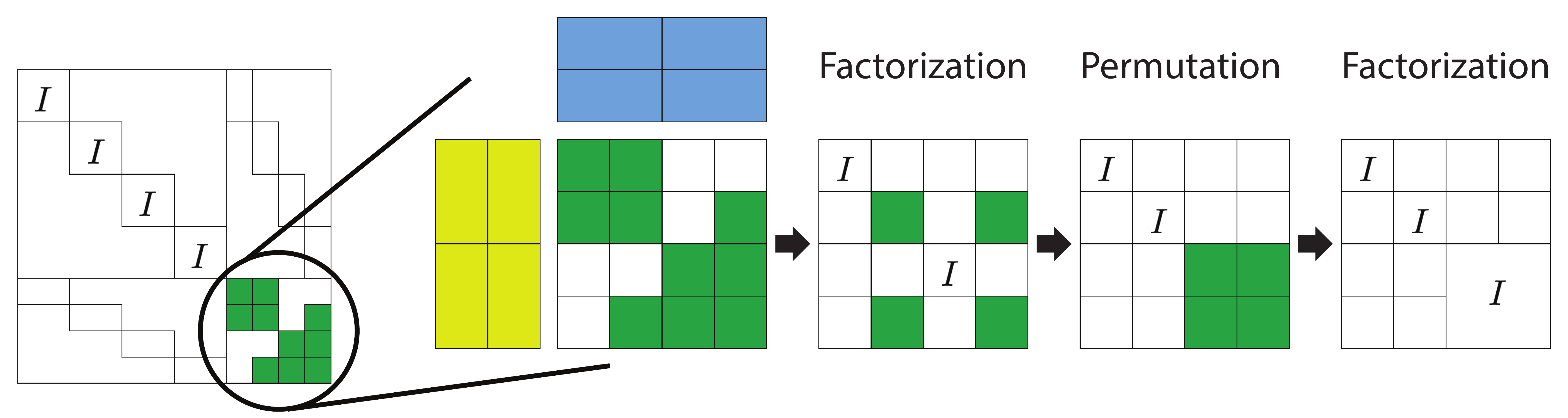}
  \caption{Upper levels of the HSS-ULV factorization.}
  \label{fig:HSS-ULV2}
\end{figure}

The HSS structure is a multi-level version of the BLR$^2$ structure with weak admissibility.
The flow of the HSS-ULV factorization is identical to that of the BLR$^2$-ULV factorization until the step in Eq. (\ref{eq:BLR2final}).
As shown in Fig. \ref{fig:HSS-ULV2}, instead of treating the leftover block as a dense block, the HSS-ULV recursively applies the same procedure to the remaining part.
The \yellowtext{column basis} and \bluetext{row basis} shown in Fig. \ref{fig:HSS-ULV2} are actually the transfer matrices \yellowtext{$U_{1;0}$} and \bluetext{$V_{1;1}$} shown in Eq. (\ref{eq:transfer}).
They can be computed from
\begin{align}
\yellowtext{[U^S_{1;0}\ U^R_{1;0}]},R =
\mathrm{QR}\left(
\yellowtext{
\begin{bmatrix}
U_{2;0} & 0 \\
0 & U_{2;1}
\end{bmatrix}^\top}
\redtext{A_{1;0,1}}
\right)\\
\bluetext{[\tp{V^S_{1;1}}\ \tp{V^R_{1;1}}]},R =
\mathrm{QR}\left(
\bluetext{
\begin{bmatrix}
V_{2;2} & 0 \\
0 & V_{2;3}
\end{bmatrix}}
\redtext{A_{1;0,1}^\top}
\right).
\end{align}
Before we compute the skeleton/redundant decomposition of the \greentext{$S$} block as in Eq. (\ref{eq:rs}), we first need to merge the \greentext{$S$} blocks in Eq. (\ref{eq:BLR2final}) as follows.
\begin{equation}
\begin{bmatrix}
S_{1;0,0} & S_{1;0,1}\\
S_{1;1,0} & S_{1;1,1}
\end{bmatrix}=
\setlength\arraycolsep{0pt}
\begin{bmatrix}
S^{SS}_{2;0,0} & S^{SS}_{2;0,1} & S^{SS}_{2;0,2} & S^{SS}_{2;0,3}\\
S^{SS}_{2;1,0} & S^{SS}_{2;1,1} & S^{SS}_{2;1,2} & S^{SS}_{2;1,3}\\
S^{SS}_{2;2,0} & S^{SS}_{2;2,1} & S^{SS}_{2;2,2} & S^{SS}_{2;2,3}\\
S^{SS}_{2;3,0} & S^{SS}_{2;3,1} & S^{SS}_{2;3,2} & S^{SS}_{2;3,3}
\end{bmatrix}
\end{equation}
Then we can decompose all the \greentext{$S$} matrices of the dense diagonal blocks as
\begin{equation}
\greentext{
\begin{bmatrix}
S^{RR}_{1;i,i} & S^{RS}_{1;i,i} \\
S^{SR}_{1;i,i} & S^{SS}_{1;i,i}
\end{bmatrix}}
=\yellowtext{\tp{[U^R_{1;i}\ U^S_{1;i}]}}
\greentext{S_{1;i,i}}
\bluetext{
\tp{
\begin{bmatrix}
V^R_{1;i} \\
V^S_{1;i}
\end{bmatrix}}}.
\end{equation}
For the low-rank blocks we have
\begin{equation}
\greentext{S^{SS}_{1;i,j}}
=\yellowtext{\tp{U^S_{1;i}}}
\greentext{S_{1;i,j}}
\bluetext{
\tp{
V^S_{1;j}}},
\end{equation}
Once all the \greentext{$S$} blocks at level 1 are computed, the same procedure shown in Eqs. (\ref{eq:lu}) to (\ref{eq:BLR2final}) is repeated.
We have shown the HSS matrix with a $4\times 4$ subdivision as an example, but the same recursive algorithm can be extended to any number of subdivision levels and any matrix size.

\subsection{$\mathcal{H}^2$-ULV factorization \rev{with dependencies}}
\begin{figure}[t]
  \centering
  \includegraphics[width=\linewidth]{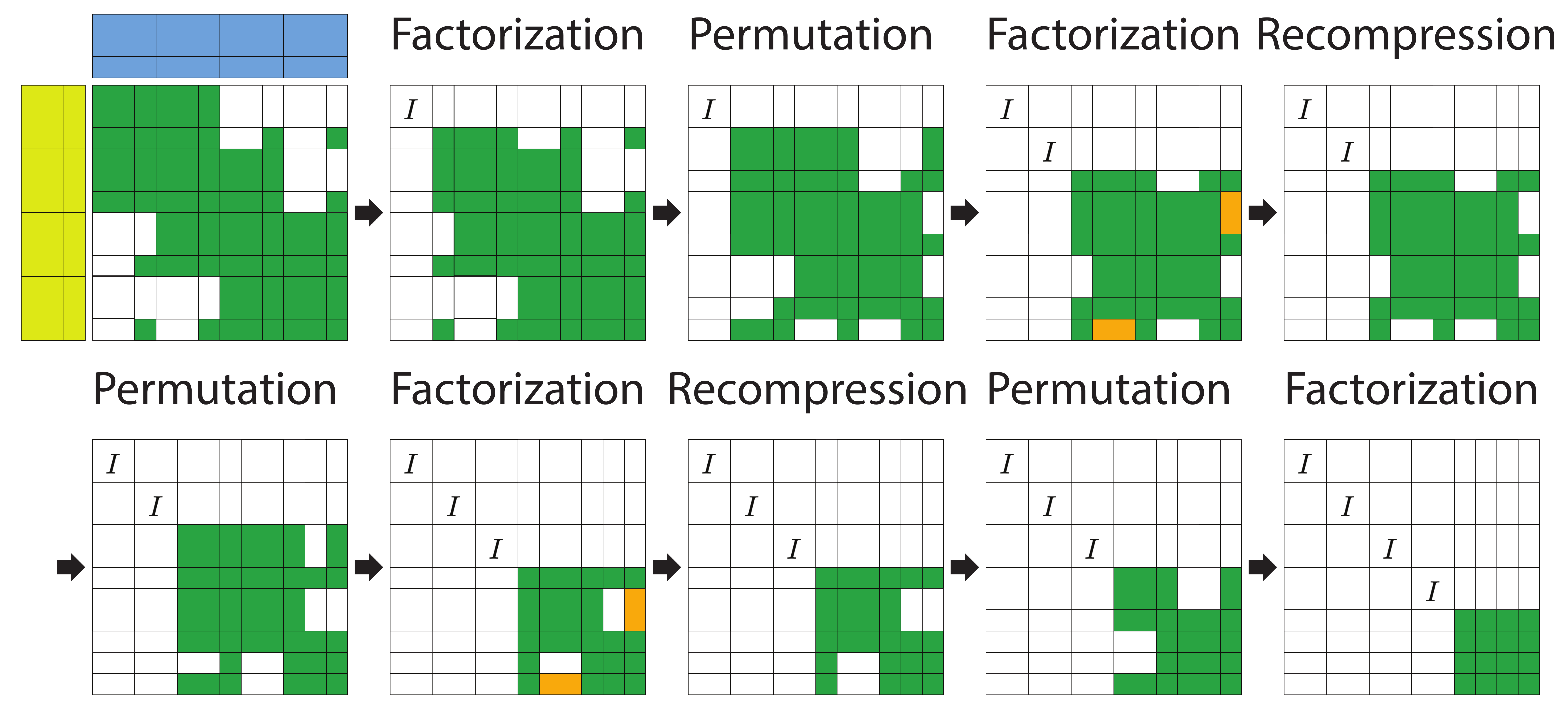}
  \caption{Flow of the $\mathcal{H}^2$-ULV factorization with dependencies.}
  \label{fig:H2-ULV}
\end{figure}

The flow of the $\mathcal{H}^2$-ULV factorization is shown in Fig. \ref{fig:H2-ULV}.
The strong admissibility of $\mathcal{H}^2$-matrices gives dense off-diagonal blocks.
\rev{The example shown in Fig. \ref{fig:H2-ULV} has a block tri-diagonal structure, but $\mathcal{H}^2$-matrices can handle any pattern of off-diagonal dense blocks.
The factorization cost will be $\mathcal{O}(N)$ as long as the number of dense blocks is $\mathcal{O}(N)$.
However, these off-diagonal dense blocks}
result in trailing sub-matrix dependencies, so the factorization at each level cannot be done in parallel like the HSS-ULV factorization.
\rev{This is one of the main reasons why HSS-ULV has been preferred over $\mathcal{H}^2$-ULV despite its suboptimal performance in 3-D problems.}
We will describe in the next section how these trailing sub-matrix dependencies can be removed, but here we will first describe the one with dependencies.
For $\mathcal{H}^2$-ULV, the stages for constructing the $\yellowtext{U}\greentext{S}\bluetext{V}$ decomposition given by Eqs. (\ref{eq:u}) to (\ref{eq:ss}), are identical to the BLR$^2$-ULV factorization. The \yellowtext{$A_{2;0,+}$} matrix in Eq. (\ref{eq:u}) looks like the concatenation of all \textit{off-diagonal} blocks in Fig. \ref{fig:notation}, but for $\mathcal{H}^2$-matrices, it is the concatenation of all \textit{low-rank} blocks.
Dense blocks are never used for the construction of the shared basis. Eq. (\ref{eq:rs}) is applied to all the dense blocks, not only the diagonal blocks.
Otherwise, the $\yellowtext{U}\greentext{S}\bluetext{V}$ decomposition is identical to that of the BLR$^2$-ULV factorization.

The factorization stage is quite different from that of the HSS-ULV, since the off-diagonal blocks create a trailing sub-matrix dependency.
This has two consequences, 1) the serialization of the factorization and permutation steps, and 2) the existence of fill-ins that result in an extra recompression step.
While the HSS-ULV shown in Fig. \ref{fig:HSS-ULV} has only a single factorization step and permutation step per level, the $\mathcal{H}^2$-ULV performs a series of factorization, recompression, and permutation steps for each block row/column.
The \orangetext{fill-in} blocks are shown in orange in Fig. \ref{fig:H2-ULV}.
These fill-in blocks can be eliminated by merging them back into the row/column basis.
For example for the fourth matrix from the left on the top row of Fig. \ref{fig:H2-ULV}, there are two \orangetext{fill-in} blocks. 
The recompression of these blocks involve the following updates to the \yellowtext{column basis} and \bluetext{row basis} as
\begin{align}
\yellowtext{[U^S_{2;2}\ U^R_{2;2}]},R &= \mathrm{QR}\left(\left[
\yellowtext{U^S_{2;2}}\greentext{S^{SS}_{2;2,0}}\hspace{10pt}
\yellowtext{
\begin{bmatrix}
U^R_{2;2}\ U^S_{2;2}
\end{bmatrix}}
\begin{bmatrix}
\orangetext{S^{RS}_{2;2,0}} \\
\greentext{S^{SS}_{2;2,0}} 
\end{bmatrix}
\right]\right) \label{eq:u2}\\
\bluetext{[\tp{V_{2;2}^S}\ \tp{V_{2;2}^R}]},R &= \mathrm{QR}\left(
\tp{
\begin{bmatrix}
\greentext{S^{SS}_{2;0,2}}\bluetext{V^S_{2;2}} \\
[\orangetext{S^{SR}_{2;0,2}}\ \greentext{S^{SS}_{2;0,2}}]
\bluetext{
\begin{bmatrix}
V_{2;2}^R \\
V_{2;2}^S
\end{bmatrix}}
\end{bmatrix}}
\right), \label{eq:v2}
\end{align}
where $\yellowtext{U^S_{2;2}}\greentext{S_{2;2,0}}$ and $\greentext{S_{2;0,2}}\bluetext{V^S_{2;2}}$ are products of the basis and the skeleton matrix.
For this particular case, there is only one low-rank block in this row/column, but in the general case it is necessary to concatenate this product for all low-rank blocks in the row/column that is being recompressed.
The operations in Eqs. (\ref{eq:u2}) and (\ref{eq:v2}) allow the \orangetext{fill-in} to be incorporated into the \yellowtext{$U$} and \bluetext{$V$} bases, so that it can be eliminated from the \greentext{$S$} matrix.
All \orangetext{fill-in} blocks are eliminated as soon as they fill-in, so the sparsity of the \greentext{skeleton} matrix is retained.
\rev{By comparing Fig. \ref{fig:HSS-ULV} and \ref{fig:H2-ULV}, it is clear that the $\mathcal{H}^2$-ULV is much more complicated and much less parallel.}

\section{$\mathcal{H}^2$-ULV without dependencies}
Existing methods without trailing sub-matrix dependencies such as GOFMM~\cite{yuDistributedMemoryHierarchicalCompression2018} and STRUMPACK~\cite{rouetDistributedmemoryPackageDense2015} \rev{are} limited to weak admissibility \rev{\textit{e.g.} HSS matrices}.
For 3-D problems, the rank of the off-diagonal blocks \rev{in an HSS matrix} grows as a function of $N$, so the $\mathcal{O}(N)$ complexity of the algorithm is lost.
On the other hand, $\mathcal{H}^2$-matrices~\cite{maDirectSolutionGeneral2019,bormVariableOrderDirectional2019} with strong admissibility do not have this problem, and are able to achieve $\mathcal{O}(N)$ even for problems with 3-D geometry. 
\rev{The $\mathcal{H}^2$-matrices having dense off-diagonal blocks is both a blessing and a curse, since this is what allows the low-rank off-diagonal blocks to have $\mathcal{O}(1)$ ranks, while it is also the source of fill-ins that ultimately serialize the factorization process.
Previous implementations of $\mathcal{H}^2$-matrix factorization such as IFMM~\cite{takahashiParallelizationInverseFast2020}} depend on coloring schemes to extract the parallelism from an inherently serial algorithm.
In the present work, we develop a method that has the best of both worlds, where we use $\mathcal{H}^2$-matrices in order to \rev{achieve linear complexity for \rev{3-D} problems}, while also \rev{removing the} dependencies \rev{on the trailing sub-matrices} during the factorization.

\begin{figure}[t]
  \centering
  \includegraphics[width=\linewidth]{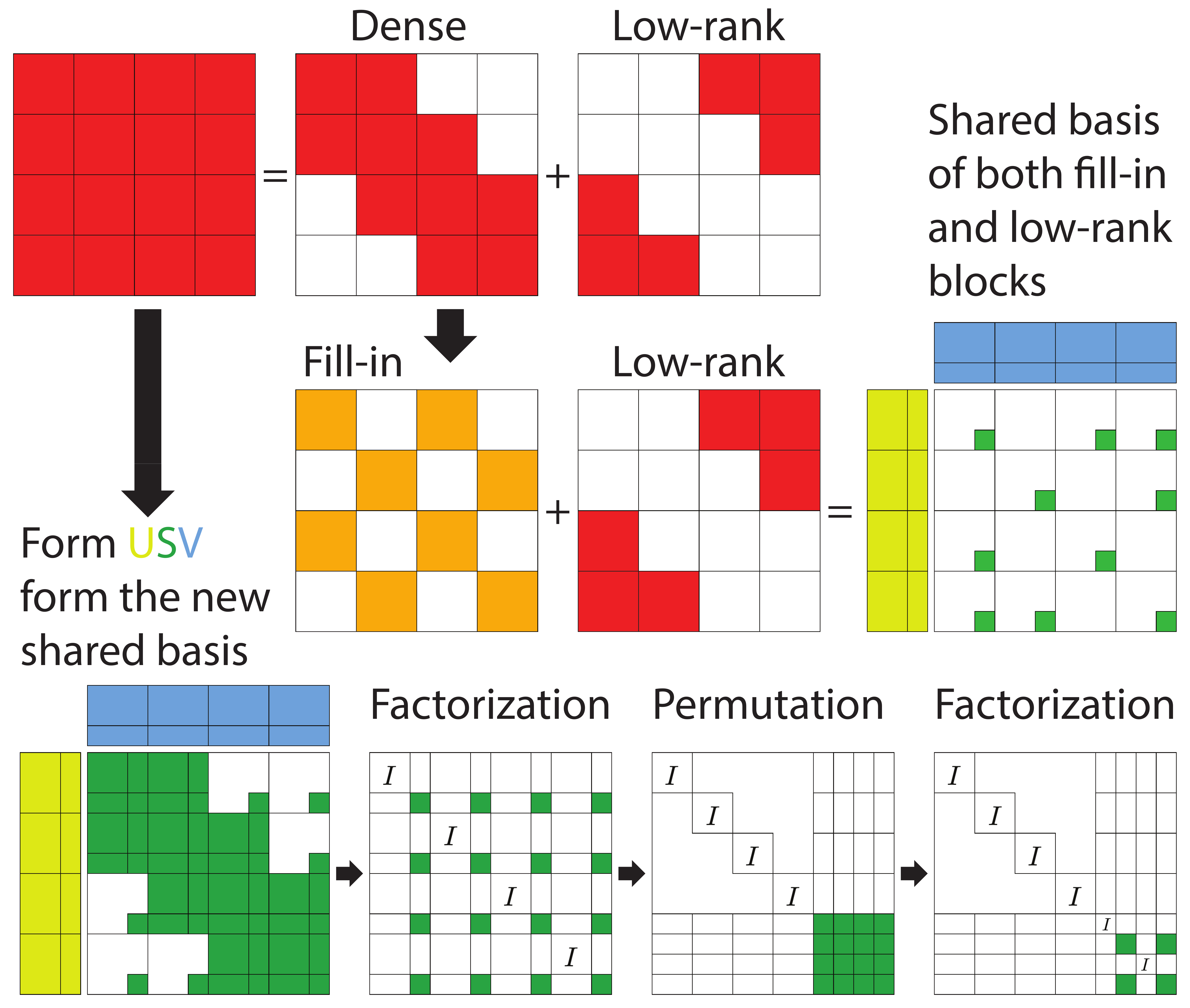}
  \caption{Flow of the $\mathcal{H}^2$-ULV factorization without dependencies.}
  \label{fig:H2-ULV2}
\end{figure}

\subsection{\rev{Fill-in blocks are low-rank}}
\rev{The nullity theorem \cite{vandebrilNoteNullityTheorem2006} states that the LU factors of a matrix have the same low-rank structure as the matrix before it is factorized.
This is only possible if all the fill-ins can be compressed back to low-rank blocks.
Existing $\mathcal{O}(N)$ dense direct solvers such as IFMM~\cite{ambikasaranInverseFastMultipole2014} and HIF~\cite{hoHierarchicalInterpolativeFactorization2016} are also based on this property of the fill-ins being low-rank.
In the present work, we also exploit this property of structured low-rank matrices to recompress the fill-in blocks to low-rank.
Furthermore, under the $\yellowtext{U}\greentext{S}\bluetext{V}$ formulation, compressing a block to low-rank is equivalent to eliminating that block for that level, as we have shown in Section \ref{sec:structured}.
Therefore, none of the blocks remain filled-in during the $\mathcal{H}^2$-ULV factorization.
They may temporarily fill-in, but are immediately eliminated through recompression to low-rank blocks.
However, it can be seen by comparing Figs. \ref{fig:HSS-ULV} and \ref{fig:H2-ULV} that the $\mathcal{H}^2$-ULV needs to sequentially factorize the \greentext{S} matrix, whereas the HSS-ULV can factorize each level in parallel.
This is due to the fact that the fill-in blocks still need to be recompressed to low-rank blocks, and the recompression requires an update to the shared basis.
Since the next block row/column cannot proceed until the shared basis is updated, this is what actually causes the serialization.
Therefore, this update to the shared basis is what prevents us from removing the dependencies on the trailing sub-matrices.
If we can somehow form a shared basis that contains all the fill-in blocks (and not just the low-rank blocks), there would be no need to update this shared basis during the factorization.
This would lead to an $\mathcal{H}^2$-ULV factorization without trailing sub-matrix dependencies.
Unlike the fill-reducing ordering in sparse direct solvers, the objective here is not to reduce the fill-ins.
It is rather to avoid the need to update the shared basis, even when some of the blocks temporarily fill-in and need to be recompressed to low-rank so that they do not remain filled-in.

\begin{figure}
    \centering
    \includegraphics[width=\linewidth]{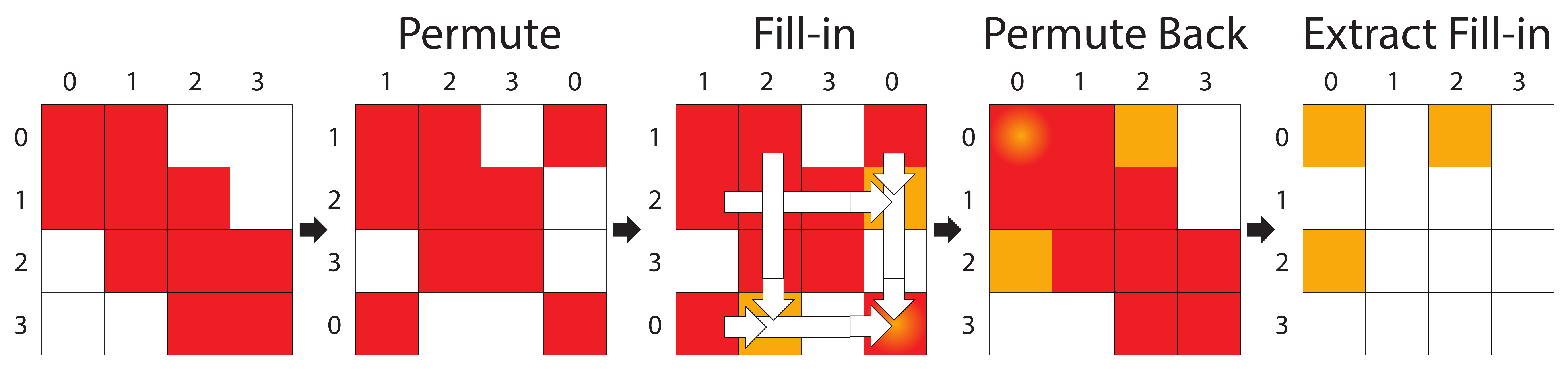}
    \caption{Fill-in of the first row and column during the pre-factorization.}
    \label{fig:fill-in}
\end{figure}

\subsection{\rev{Pre-computing the fill-ins per block row/column}}
Our method consists of two separate factorization phases -- one for dense blocks in the original matrix $\redtext{A}$ to compute the fill-ins, and another for the $\greentext{S}$ matrix after the shared basis has been constructed. 
Both phases can be done without dependency on trailing sub-matrices.
The flow of computation is shown in Fig. \ref{fig:H2-ULV2}.
We first decompose the matrix into the dense blocks and low-rank blocks.
Then, the fill-ins that occur during the factorization of the dense blocks is computed using the method shown in Fig. \ref{fig:fill-in}.
In order to compute all potential fill-in blocks in a given block row/column, we permute that block row to the bottom and block column to the far right.
We then compute all possible \orangetext{fill-ins} from all other rows/columns into that block row/column.
In this particular case, the $A_{2;0,2}$, $A_{2;2,0}$, and $A_{2;0,0}$ blocks will fill-in.
We do not accumulate the fill-in blocks into the dense blocks during this process, but store them separately.
We then permute the block row/column back to their original position along with the shared basis.
Actually, the permutation is done only for the sake of explaining the fill-in process in an intuitive way, and this process can be implemented without the permutation if one can identify which blocks will fill-in without permuting them.
This process can be executed in parallel for all block rows/columns, since they do not depend on each other.
Note that the factorization of the diagonal blocks and the triangular solves of the dense off-diagonal blocks are not redundantly computed for the same block more than once during this process.
The resulting \orangetext{fill-in} matrix is shown at the center of Fig. \ref{fig:H2-ULV2}.}

\subsection{\rev{A shared basis for both fill-in and low-rank blocks}}
\rev{The next step is to form a shared basis between the fill-in and low-rank blocks for each block row/column.
This shared basis will contain all the information necessary to compress the low-rank blocks initially, and any of the fill-in blocks that arise during the factorization.
All fill-in and low-rank matrices for a given row/column are concatenated to form the shared row/column bases
\begin{align}
\yellowtext{[U^S_{2;i}\ U^R_{2;i}]},R &= \mathrm{QR}([\orangetext{F_{2;i,+}\ \redtext{A_{2;i,+}}}]) \label{eq:fillu}\\
\bluetext{[\tp{V_{2;j}^S}\ \tp{V_{2;j}^R}]},R &= \mathrm{QR}([\orangetext{\tp{F_{2;+,j}}}\ \redtext{\tp{A_{2;+,j}}}]),\label{eq:fillv}
\end{align}
where $\orangetext{F_{2;i,+}}$ and $\redtext{A_{2;i,+}}$ are the concatenation of all \orangetext{fill-in} blocks and \redtext{low-rank} blocks in the $i$th row, respectively.
Similarly, $\orangetext{F_{2;+,j}}$ and $\redtext{A_{2;+,j}}$ are the concatenation of all \orangetext{fill-in} blocks and \redtext{low-rank} blocks in the $j$th column, respectively.
The center row in Fig. \ref{fig:H2-ULV2} corresponds to the operation in Eqs. (\ref{eq:fillu}) and (\ref{eq:fillv}).
The procedure is identical to the construction of the shared basis for the low-rank blocks besides the inclusion of the fill-in blocks.
These bases are then used to compute the \greentext{$S$} matrices for the dense blocks using Eq. (\ref{eq:rs}) and the low-rank blocks using Eq. (\ref{eq:ss}).}
\rev{Once the $\yellowtext{U}\greentext{S}\bluetext{V}$ decomposition is formed for the shared basis that incorporates both the fill-in and low-rank blocks, the factorization of the $\mathcal{H}^2$-ULV can proceed just as the HSS-ULV, as shown in Fig. \ref{fig:H2-ULV2}.}
The same procedure shown in Eqs. (\ref{eq:lu}) to (\ref{eq:BLR2final}) can be applied here as well.
After the leaf level has been processed, the matrix can be permuted to cluster the remaining skeleton parts, which can be solved recursively by applying the same procedure at each level.

In summary, the key ideas that make it possible to remove the trailing sub-matrix dependency even for strong admissibility are:
\rev{
\begin{itemize}
    \item Fill-in blocks are always low-rank. This has been demonstrated in previous studies~\cite{ambikasaranInverseFastMultipole2014,hoHierarchicalInterpolativeFactorization2016}, and is not a unique claim of this paper.
    \item The fill-in blocks can be pre-computed and shared bases can be formed to incorporate both the fill-in and low-rank blocks. To the extent of our knowledge, this has not been done before.
    \item With this new shared basis that incorporates the fill-in blocks, there is no need to update the shared basis during the factorization of $\mathcal{H}^2$-matrices, which allows us to remove the trailing sub-matrix dependency.
\end{itemize}
}

\subsection{\rev{Parallel implementation on distributed memory}}
\rev{
Previous attempts to parallelize $\mathcal{H}^2$-ULV required coloring~\cite{takahashiParallelizationInverseFast2020} or the use of task-based runtime systems~\cite{cambierTaskbasedDistributedParallel2021} in order to extract parallelism from an inherently serial algorithm.
The algorithm is inherently serial because there are trailing sub-matrix dependencies, where the blocks in the lower-right need to wait for the blocks in the upper-left to finish computing the GETRF (LU factorization), TRSM (traiangular solve), and GEMM (matrix multiplication for the Schur complements).
Our proposed method removes this dependency by pre-computing a shared basis that does not need to be updated during the factorization.
The only dependency that remains is the one between the levels.
(By "level", we mean the level of granularity in the hierarchically sub-divided matrix structure.)
We can somewhat relax this level-wise dependency as well, by computing the upper levels redundantly on multiple processes.

\begin{figure}
    \centering
    \includegraphics[width=\linewidth]{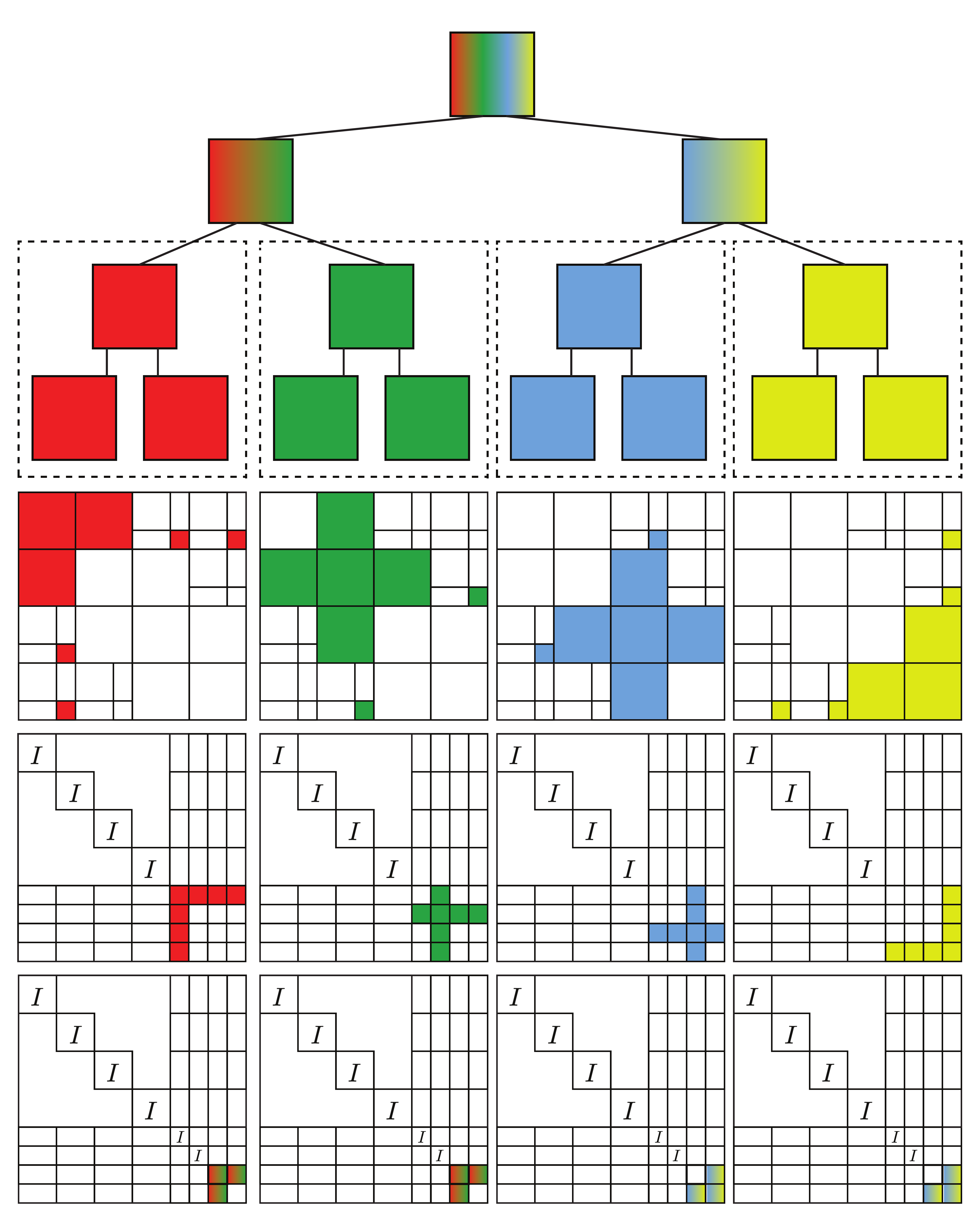}
    \caption{Partitioning of the $\mathcal{H}^2$-matrix. Upper levels are computed redundantly by multiple processes.}
    \label{fig:h2-dist}
\end{figure}

The resulting partitioning scheme is shown in Fig. \ref{fig:h2-dist}, where the different colors represent the different partitions.
In this example, each partition has only one block row/column, but there could be many block rows/columns in each partition.
The cross shape of the partitions comes from partitioning the underlying geometry, where a certain block row and block column correspond to a sub-domain in the geometry.
At the upper levels in the matrix the colors are mixed, which represents the redundant storage and computation of these blocks.
This redundancy is the key to achieving good scalability at the upper levels, because it utilizes the processes that will be idle otherwise, while eliminating the need to communicate the results to each other.

Since it is always possible to split the range of processes in half (for odd numbers roughly half), the process tree shown in Fig. \ref{fig:h2-dist} is always a full binary tree, regardless of the underlying geometry or the type of matrix.
The rows and columns of the $\mathcal{H}^2$-matrix also form a full binary tree, which is usually deeper than the process tree.
This means that the lower levels of the row/column tree are grafted to the leaves of the process tree as shown in Fig. \ref{fig:h2-dist}.
When the factorization reaches the level where more than one process owns the block, the two child blocks exchange their information through an the \texttt{Allgather} collective with a split communicator, which exchanges information between the pair of child blocks that are being merged.
At even higher levels in the process tree we can also exchange the necessary information between more processes through an \texttt{Allgather} collective with a communicator that is split accordingly.
Note that for these \texttt{Allgather} collectives with split communicators, each process only communicates with one other process at any given level.
We are essentially performing a tree \texttt{Allgather} through a hierarchy of communicators that are split according to the process tree.
However, the size of the blocks that are gathered become smaller as the level increases, since $3/4$ of the blocks get eliminated at each level for the upper levels of the $\mathcal{H}^2$-matrix.
}

\begin{figure}    
    \centering
    \begin{subfigure}{\linewidth}
        \includegraphics[width=\textwidth]{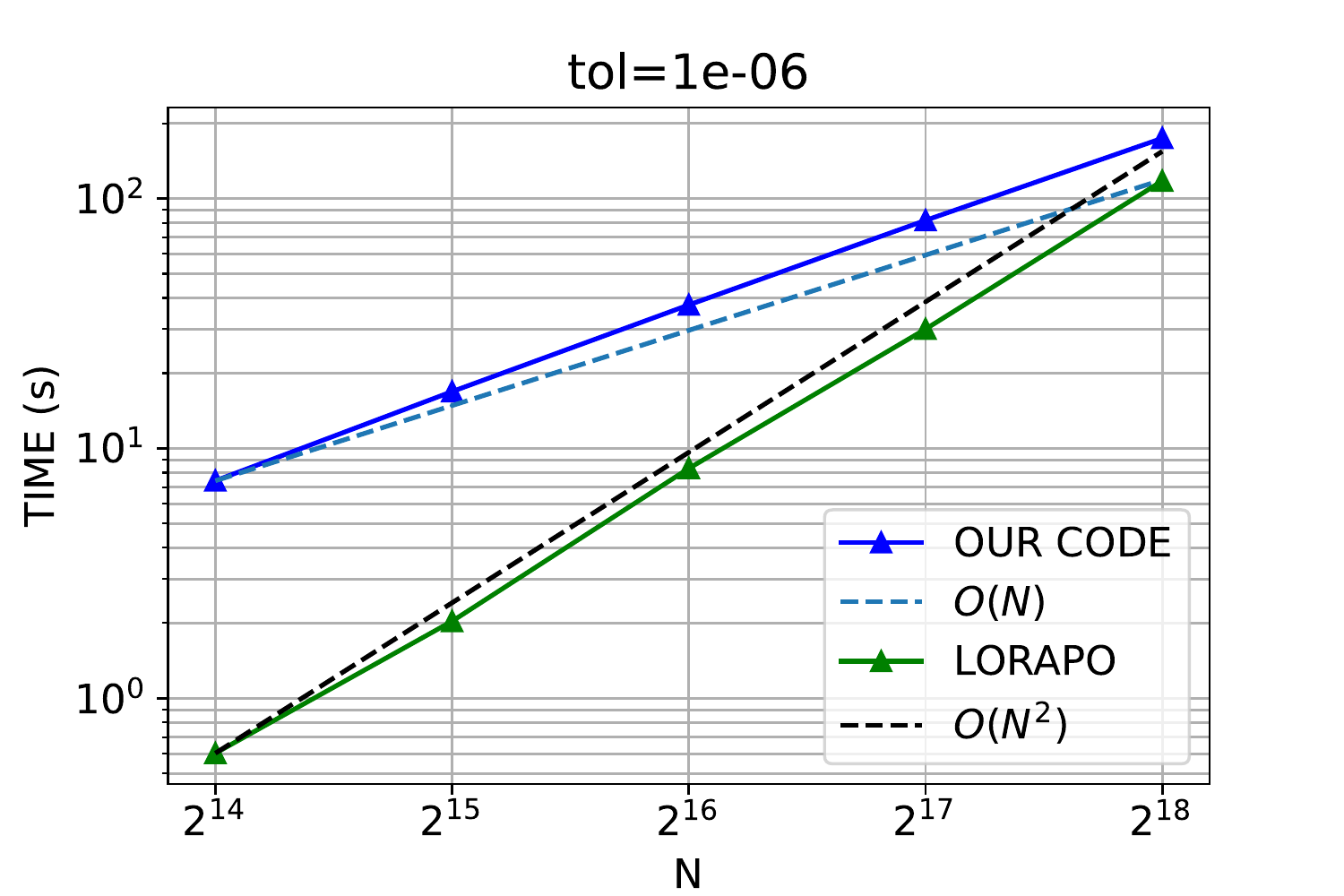}
        \caption{Relative error of $10^{-6}$.}
        \label{fig:single-core-1e-6}
    \end{subfigure}
    \begin{subfigure}{\linewidth}
        \includegraphics[width=\textwidth]{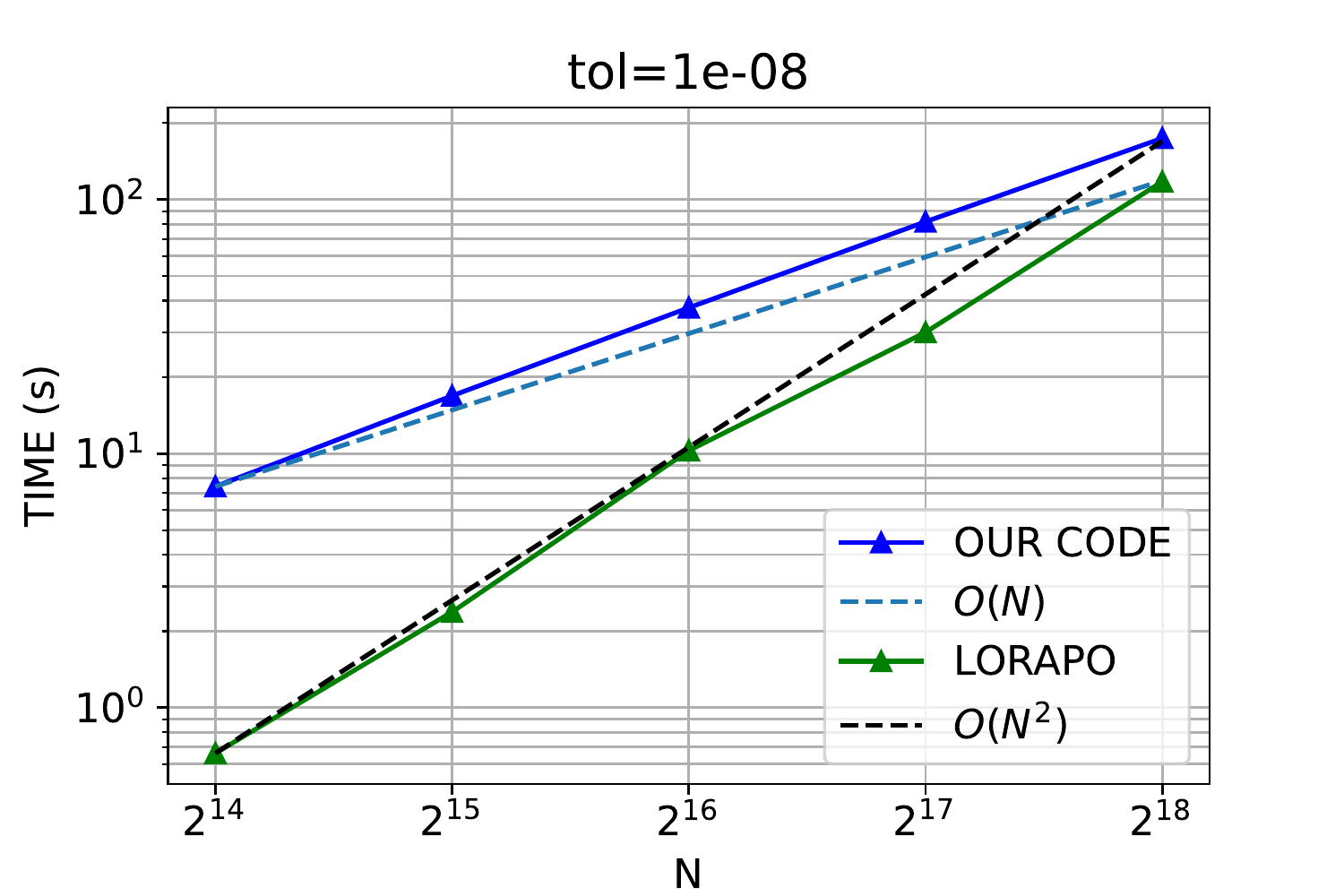}
        \caption{Relative error of $10^{-8}$.}
        \label{fig:single-core-1e-8}
    \end{subfigure}
    \caption{Comparison of LORAPO vs. our code on a single core for different \rev{problem sizes and} relative error.}
    \label{fig:single-core-comparison-lorapo-vs-our-code}
\end{figure}

\rev{
\section{Tests on a simple geometry with Laplace potential}
We first test our algorithm on a simple geometry with particles uniformly distributed inside a 3-D unit cube.
We assume unit charges on each particle.
We use the Green's function solution of the Laplace equation as our kernel.
\begin{equation}
\Phi_i=\frac{q_j}{4\pi r_{ij}},
\end{equation}
where $r_{ij}$ is the Euclidean distance between the points $\mathbf{x}_i$ and $\mathbf{x}_j$, and $q_j$ is the charge.
The Green's function matrix $G_{ij}=\frac{1}{4\pi r_{ij}}$ results in a dense but rank-structured matrix.
}

\subsection{Experimental setup}
We perform tests on a single node with 2xAMD EPYC 7742 CPUs, each with 64 physical cores running at 3.3 GHz.
The total available memory of the node is 1000GB, evenly split between both CPUs.
We use GCC 9.3.0 as our compiler and openMPI 4.0.3 for the MPI processes, and link to Intel MKL 2020.1.

We compare our implementation with LORAPO \cite{caoExtremeScaleTaskBasedCholesky2020}, an  adaptive-rank BLR Cholesky factorization using the PaRSEC PTG \cite{caoLeveragingPaRSECRuntime2021} runtime system for achieving asynchronous parallelism.
Although LORAPO uses mixed precision for the low rank parts \cite{abdulahAcceleratingGeostatisticalModeling2021}, we switch it off for these these tests and use only double precision for all our experiments. 
We found the most optimal way to run our code was by attaching a single process to each physical core, whereas LORAPO is run by spawning one process per node and attaching threads to each physical core.
We report our strong scaling experiments by only reporting the number of cores being utilized.
In all of our results, the relative L2 error is calculated by comparing the accuracy of the solution obtained using our method to the one obtained using a dense LU factorization from LAPACK.

Fig. \ref{fig:single-core-comparison-lorapo-vs-our-code} shows time taken for factorization for our code vs. LORAPO on a single core with varying problem sizes.
Our algorithm is able to scale linearly (ideal $O(N)$ scaling is shown by the blue dotted line) for all problem sizes.
It can be seen that although our code scales almost linearly, LORAPO has a better time to solution for most problem cases in spite of having a time complexity of $O(N^2)$, which is denoted by the black dotted line.

\begin{figure}
    \centering
    \includegraphics[width=\linewidth]{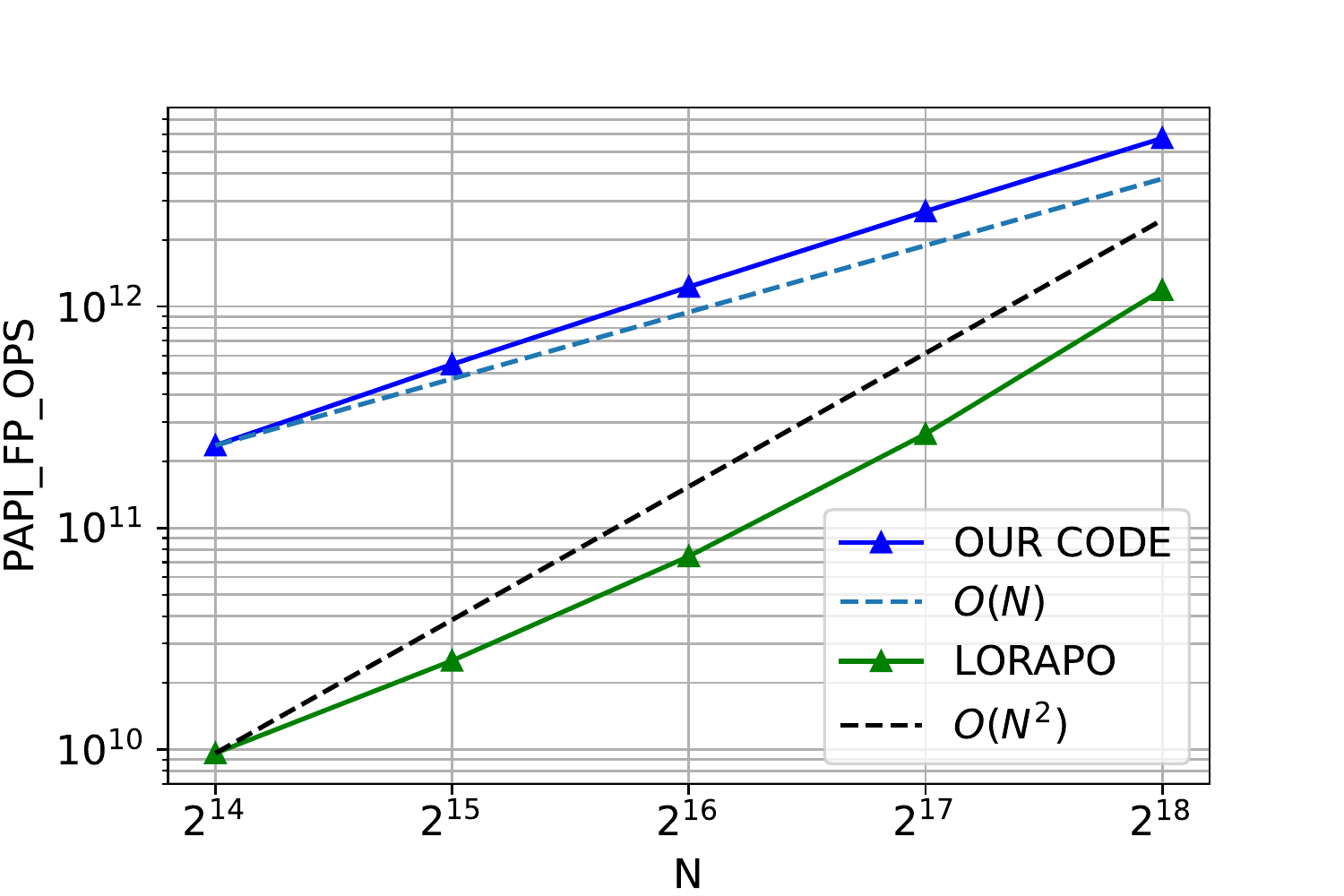}
    \caption{Comparison of PAPI\_FP\_OPS between our code and LORAPO for an 
    accuracy of $10^{-8}$ and the same tile sizes as used in Fig. \ref{fig:single-core-1e-8}}
    \label{fig:papi-fp-ops-comparison}
\end{figure}

Proof that the faster performance of LORAPO is purely due to the algorithmic
computations and not due to some anomaly in our code can be seen in 
Fig. \ref{fig:papi-fp-ops-comparison}, where we compare the number of floating point
operations performed for the same single core experiment shown in Fig. \ref{fig:single-core-1e-8}.
\rev{ULV based factorization requires more flops compared to a normal factorization.
Applying the \yellowtext{U} and \bluetext{V} basis to the dense blocks of the matrix has a cost similar to that of factorizing these dense blocks, which results in a large extra cost for ULV.
In addition, sharing and nesting bases generally results in a larger rank than what the individual low-rank blocks have.
BLR takes advantage of being able to independently compress each low-rank block, so that their rank can be minimized to save flops.
In upper levels of the $\mathcal{H}^2$-matrix, we have reported seeing a rank that is as high as 180, whereas BLR uses a maximum of rank 50 at the leaf.
A combination of these factors is what leads to the greater number of flops for our method under serial and small problem size settings.
}

\begin{figure}
    \centering
    \begin{subfigure}{\linewidth}
        \includegraphics[width=\linewidth]{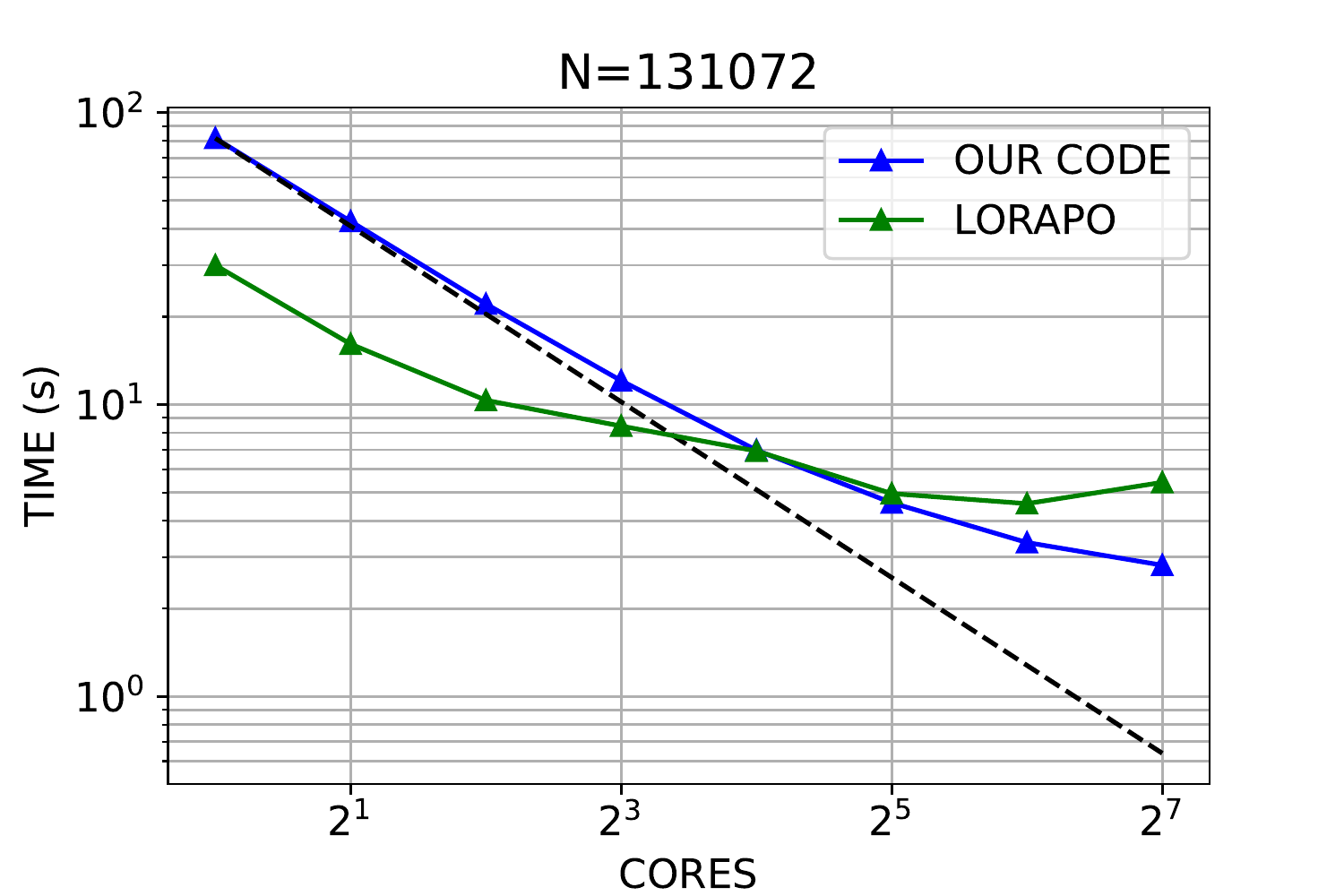}
        \caption{N=131072}
        \label{fig:strong-scaling-131072}
    \end{subfigure}
    \begin{subfigure}{\linewidth}
        \includegraphics[width=\linewidth]{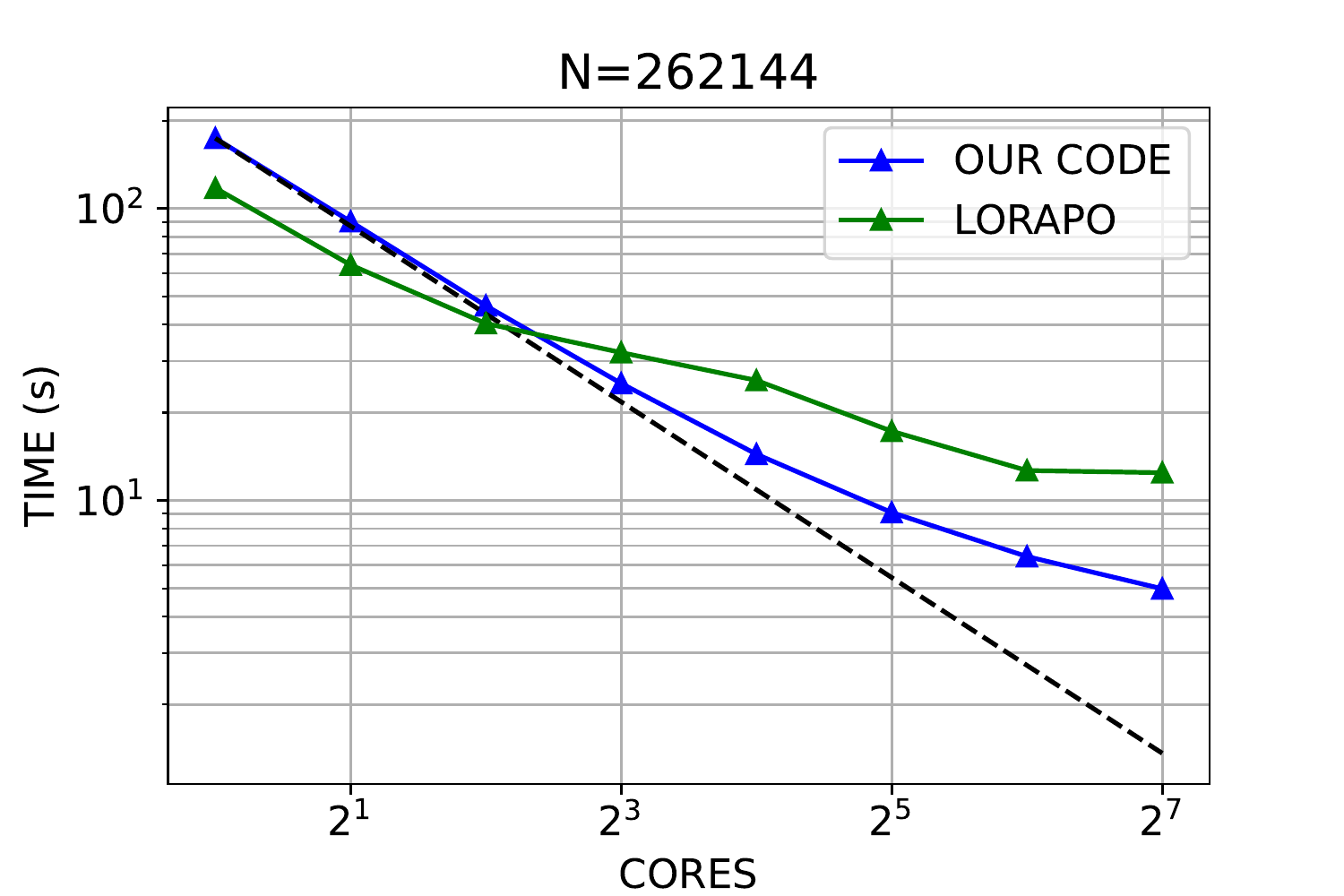}
        \caption{N=262144}
        \label{fig:strong-scaling-262144}
    \end{subfigure}
    \caption{Strong scaling experiments for various problem sizes on a single node utilizing
    upto 128 cores. The dotted black line shows perfect scaling.}
    \label{fig:strong-scaling-single-socket}
\end{figure}

\subsection{Strong scalability of shared memory parallelism}
Fig. \ref{fig:strong-scaling-single-socket} shows the strong scaling of our code vs. LORAPO for
a constant problem size of $N=131072$ in Fig. \ref{fig:strong-scaling-131072} and for $N=262144$
in Fig. \ref{fig:strong-scaling-262144}. The accuracy for all experiments is constant at $10^{-8}$.
The advantage of the inherent parallelism of our algorithm 
described in the above sections can be seen here. Our algorithm is able to beat LORAPO when using
a large number of cores, even though LORAPO uses a runtime system for asynchronous parallelism,
and uses fewer floating-point operations.

\begin{figure}
    \centering
    \includegraphics[width=\linewidth]{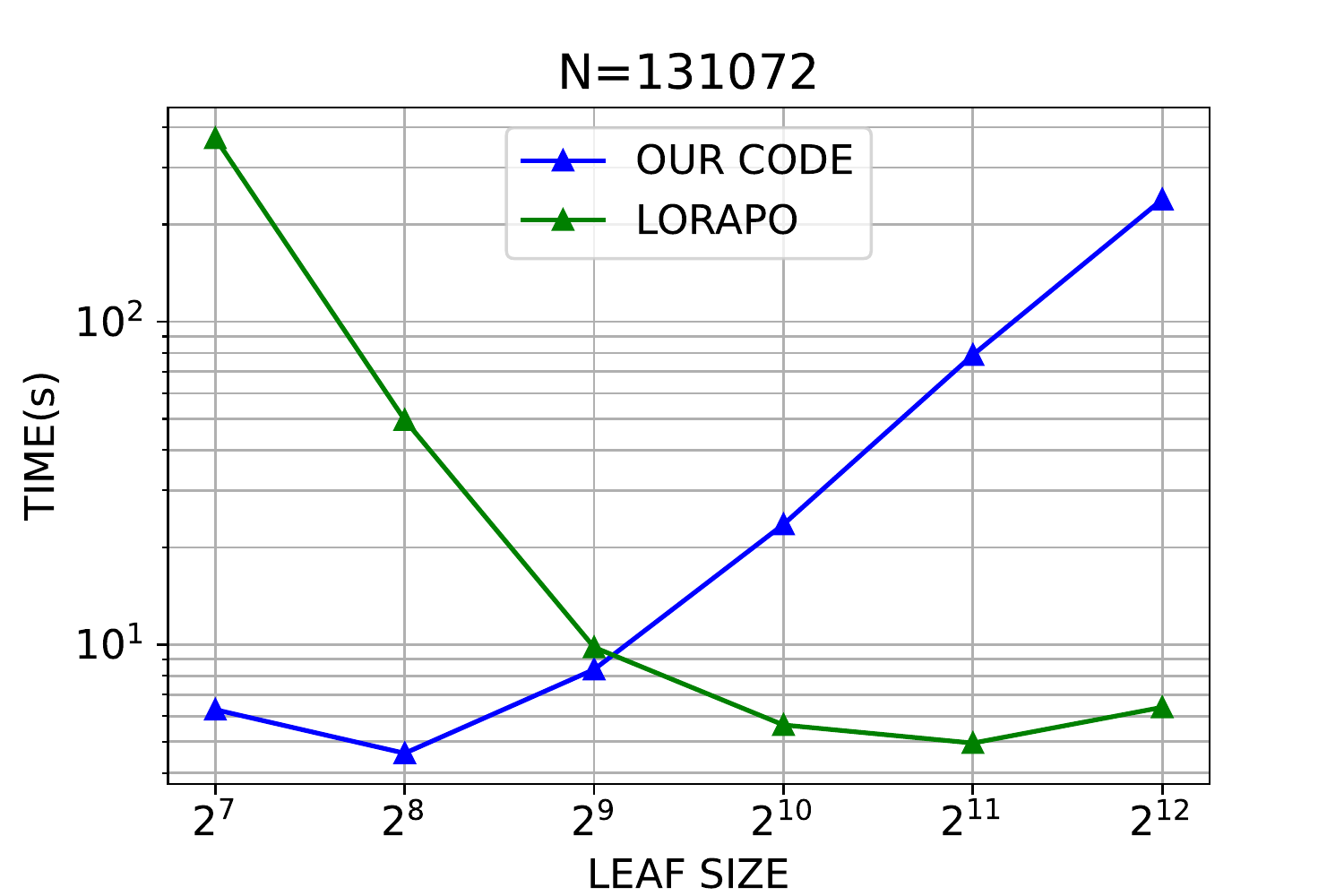}
    \caption{Impact of change the leaf size for LORAPO and our code for the same problem size ($N=131072$ and available resources (32 cores).}
    \label{fig:tile-size-impact}
\end{figure}

The reason behind poor scaling of LORAPO in Fig. \ref{fig:strong-scaling-131072} 
can be seen from the trace of the computation
taken for 64 cores on a problem size of $N=131072$ with a leaf size of $1024$. The red tasks in the trace are 
overhead introduced by the run time system PaRSEC , and the green tasks are actual useful
computation, i.e. the time actually spent in executing the tasks. The chief reason behind the
poor strong scaling is that the sizes of the tasks are too tiny to overcome the overhead 
of PaRSEC. The sizes of the red tasks are almost similar to the sizes of the useful computation.
As a result of the dependencies introduced by the Cholesky factorization algorithm, new tasks
do not become available fast enough, thus leading to poor performance.

\begin{figure}
    \centering
    \includegraphics[width=\linewidth]{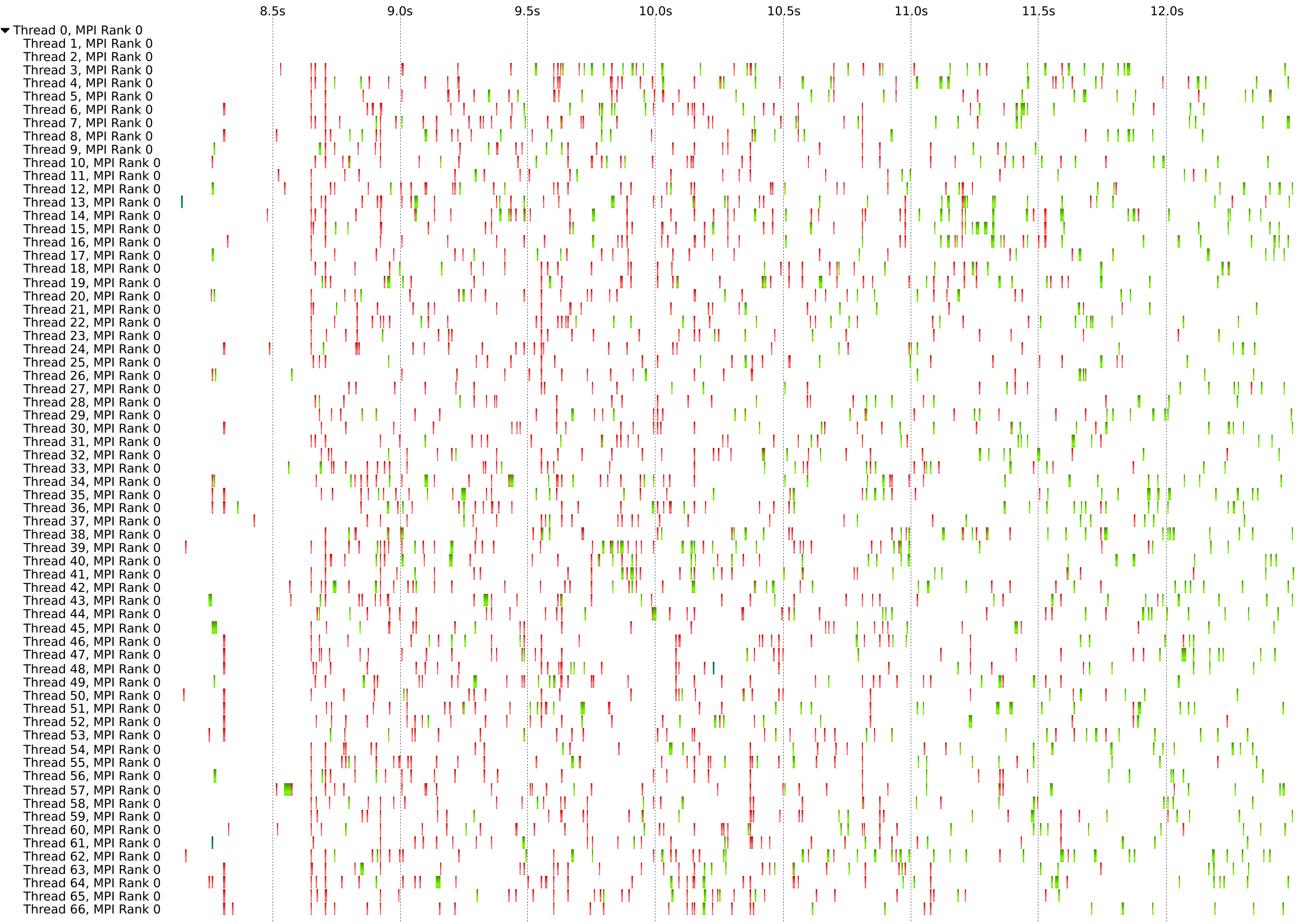}
    \caption{Trace visualization of LORAPO for a problem size of $131072$ using 64 physical cores.
    The first two threads are reserved by PaRSEC for monitoring, task submission and communication purposes
    leading to a total of 66 threads showing in the trace. The red tasks are run time system
    overhead and the green tasks are useful computation.}
    \label{fig:trace-64-cores-131072}
\end{figure}

\subsection{Variation of the leaf size}
The leaf size for LORAPO is the size of each block in the block low rank matrix.
For our code it is the number of particles present in the leaf node of the tree. Varying these
parameters, while keepin the number of cores constant at $32$ and the problem
size constant at $N=131072$ leads to changes in the time to solution as shown in 
Fig. \ref{fig:tile-size-impact}. It can be seen that LORAPO reaches an optimal
execution time as the leaf size is increased until $2048$, whereas our code
is most optimal when the leaf size is $256$.

Our algorithm uses a tree structure for representing the $\mathcal{H}^2$-matrix, and
increasing the leaf size decreases the height of the tree, which increases
the amount of computation that has to be performed for the factorization.
The increased leaf size also increases the amount of work performed per process
and reduces the available parallelism, thus leading to an increase in run time
as the leaf size is increased.
The impact of leaf size is exactly the reverse for LORAPO due to considerations
of optimal computation for the block low rank matrix structure, adaptive rank
capability of the factorization, run-time system overhead and available
parallelism. The run time increases slightly as LORAPO approaches tile size
$4096$ due to reduction of available parallelism.

\begin{figure}[t]
  \centering
  \includegraphics[width=\linewidth]{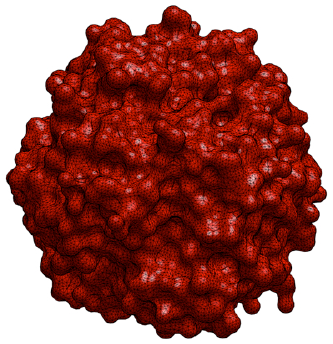}
  \caption{Boundary element mesh on a single hemoglobin.}
  \label{fig:bem}
\end{figure}

\begin{figure}
    \centering
    \includegraphics[width=\linewidth]{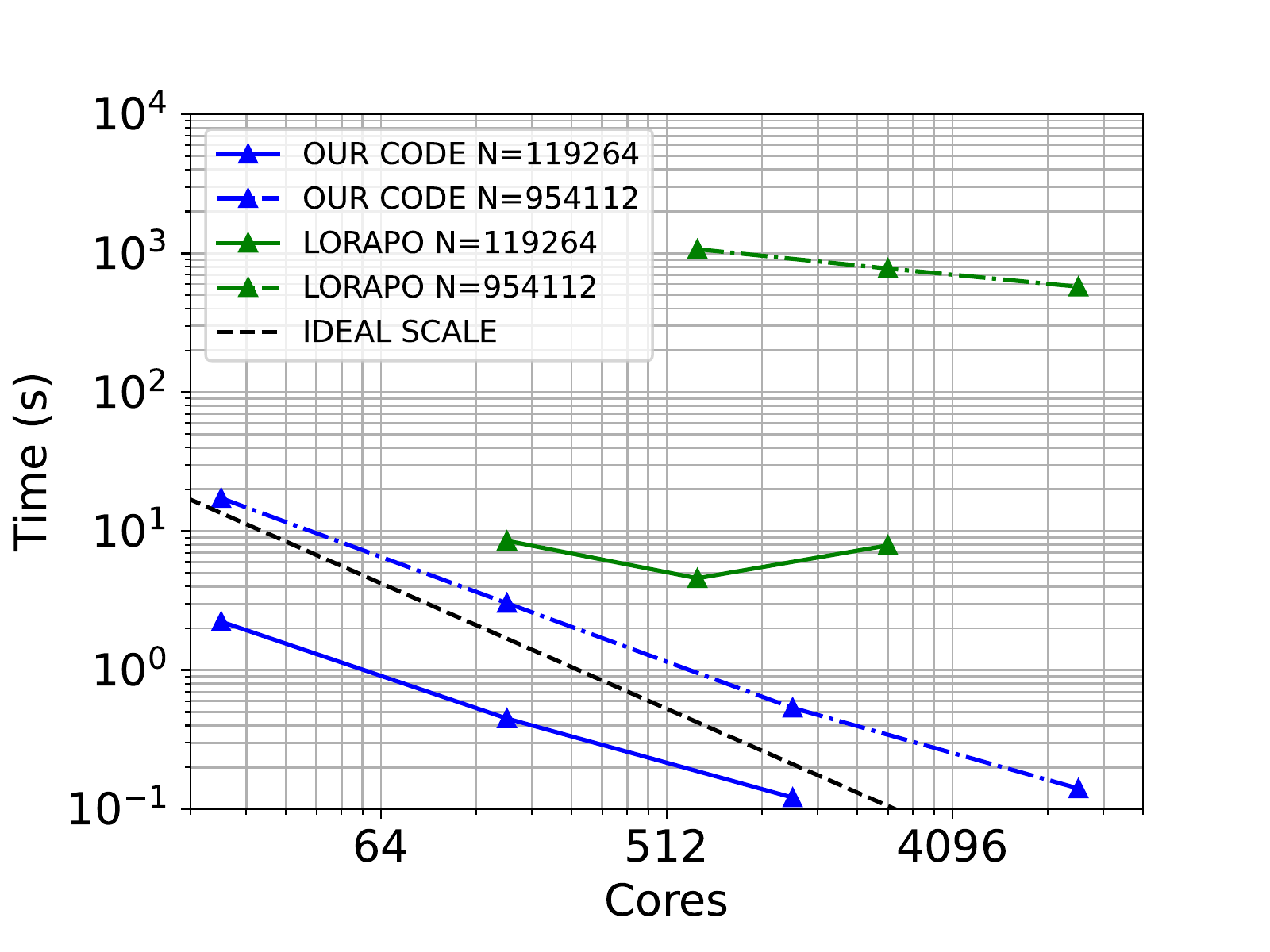}
    \caption{Strong scaling experiments on multiple nodes for different problem sizes of the hemoglobin boundary element problem using our code and LORAPO.}
    \label{fig:fac_time}
\end{figure}

\rev{
\section{Tests on complex geometry with Yukawa potential}
We now conduct experiments using an actual boundary element application in implicit solvent bio-molecular electrostatics. 
This is different from classical molecular dynamics where the water molecules are computed explicitly, adding millions more degrees of freedom to the simulation.
In implicit solvent methods, the molecules and the solvent are treated as continuous dielectric media with different dielectric constants.
This jump in the dielectric constant across the surface of the molecule is accounted for by placing a boundary element mesh over the surface.
An example of the boundary element mesh we use in the current computations is shown in Fig. \ref{fig:bem}.
For even larger problem sizes, we use a crowded environment of many hemoglobin as shown in Fig. \ref{fig:array}.
We use a collocation boundary element method, which essentially turns this mesh into a cloud of points.
We use a 3-D k-means clustering to partition those cloud of points to form the leaf blocks of the $\mathcal{H}^2$-matrix.
The flexibility of k-means clustering allows us to enforce the number of clusters to always be a power of two.
We found that this works much better than space-filling curves for partitioning points on the surface of a complex geometry.

The potential we use here is the Yukawa potential, also known as the screened Coulomb potential, which takes the following form.
\begin{equation}
\Phi_i=\frac{q_iq_j}{4\pi\varepsilon_0 r_{ij}}\exp(-\alpha m r_{ij}),
\end{equation}
where $r_{ij}$ is the Euclidean distance between the points $\mathbf{x}_i$ and $\mathbf{x}_j$, and $q_i$, $q_j$ are the charges at those points, $\alpha$ is a scaling constant, $m$ is the mass of the particle, and $\varepsilon_0$ is the permittivity.

\subsection{Experimental setup}
The distributed memory tests are performed on the ABCI machine at AIST, with a total of 1088 compute nodes.
Each node is configured with 2xIntel Xeon Gold 6148 CPUs, each with 20 physical cores running at 2.2 GHz each.
The total available memory per node is 384GB, which is evenly split between the two CPUs. We use GCC 11.2.0 compiler with Intel MPI 2021.5 and link to Intel MKL 2022.0.

\begin{figure}[t]
  \centering
  \includegraphics[width=0.95\linewidth]{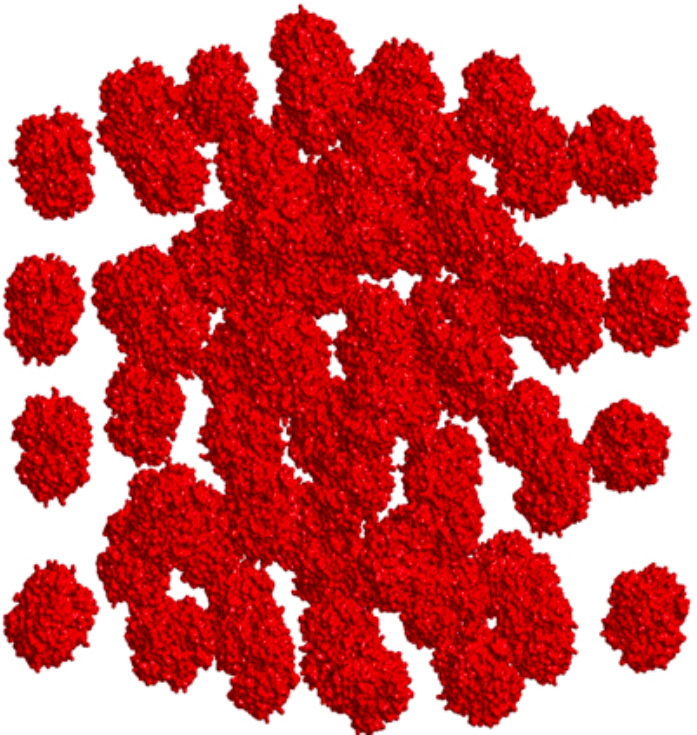}
  \caption{A crowded environment of 64 hemoglobin.}
  \label{fig:array}
\end{figure}

\subsection{Strong scalability of the distributed memory parallelism}
The results of our strong scalability tests on up to 10,240 cores are shown in Fig. \ref{fig:fac_time}.
The coverage of data points is limited for various reasons, where our code ran out of memory for some cases, while LORAPO crashed for others.
The $N=119,264$ case has a size similar to the previous experiments on the simple geometry, and on 160 cores our code is already faster.
This agrees with our results in the previous section.
When the problem size is increased to $N=954,112$, the difference between our code and LORAPO becomes much larger.
This is due to the difference between the linear complexity of our code vs. the quadratic complexity of LORAPO.
$N$ is increased roughly an order of magnitude between these two experiments, so the factorization time of our code increases and order of magnitude, while the time of LORAPO increases two orders of magnitude.
One can see that our code scales better than LORAPO as well, where the multi-node scalability of LORAPO is rather poor, while our code continues to scale well up to 10,240 cores.
If we compare the factorization time for the $N=954,112$ case on 10,240 cores between LORAPO and our code, we see that our code is approximately 4,700 times faster.
}

\section{Conclusion}
We have developed an algorithm that removes the trailing sub-matrix dependencies from $\mathcal{H}^2$-ULV factorization.
This results in an algorithm that can factorize a dense matrix in $\mathcal{O}(N)$ time that is highly parallel.
Existing work such as GOFMM and STRUMPACK use the ULV factorization to remove the trailing sub-matrix dependencies, but are limited to the HSS structure with weak admissibility.
When HSS matrices are applied to problems with 3-D geometry, \rev{the rank of the off-diagonal blocks grows as a function of $N$, which leads to suboptimal complexity.}
Our method is based on the $\mathcal{H}^2$-matrix with strong admissibility, \rev{which can handle 3-D problems in $\mathcal{O}(N)$ time.
One disadvantage of $\mathcal{H}^2$-matrices is the fill-ins, which do not exist in HSS matrices.
There are existing methods that recompress these fill-ins so that the resulting LU factors have the same low-rank structure.
However, this recompression results in an update to the shared basis, and introduces a dependency for trailing sub-matrices.
Previous methods to factorize $\mathcal{H}^2$-matrices were not able to remove this dependency on trailing sub-matrices.

The present work is able to remove this dependency on trailing sub-matrices by pre-computing all the fill-ins and including them in the shared basis before forming the $\yellowtext{U}\greentext{S}\bluetext{V}$ decomposition.
This results in a scalable dense direct solver that can handle 3-D problems in linear time without trailing sub-matrix dependencies.}
\rev{We compared our results with a block low-rank solver LORAPO, which does have trailing sub-matrix dependencies, but uses PaRSEC to alleviate the dependency issue.
Our experiments on a single node showed that our method does much more operations compared to LORAPO when the problem size is small, but our method has $\mathcal{O}(N)$ complexity, whereas LORAPO has $\mathcal{O}(N^2)$ complexity.
Therefore, our method becomes faster as the problem size grows.
Experiments on multiple nodes using up to 10,240 cores showed that our method indeed becomes much faster than LORAPO in this regime, especially for larger problem sizes.
For the problem size of close to a million and on 10,240 cores, our method showed a speed up of 4,700 fold over LORAPO.
}


\bibliographystyle{IEEEtran}
\bibliography{rioyokotalab}

\end{document}